\documentclass[12pt]{amsart}
\usepackage{graphicx}
\newtheorem{thm}{Theorem}[section]
\newtheorem{cor}[thm]{Corollary}
\newtheorem{lem}[thm]{Lemma}
\newtheorem{prop}[thm]{Proposition}
\theoremstyle{definition}

\theoremstyle{remark}
\newtheorem{rem}[thm]{Remark}
\numberwithin{equation}{section}

\begin{document}

\title[Balanced Metrics and Chow Stability]{Balanced Metrics and Chow Stability of Projective Bundles over K\"ahler Manifolds}%
\author{REZA SEYYEDALI}%
\address{Johns Hopkins University, Department of Mathematics}%
\email{seyyedali@math.jhu.edu}%

\thanks{}%
\subjclass{}%
\keywords{}%

\date{January 1, 2009}
\begin{abstract}

In 1980, I. Morrison proved that slope stability of a vector
bundle of rank $2$ over a compact Riemann surface implies Chow
stability of the projectivization  of the bundle with respect to
certain polarizations. Using the notion of balanced metrics and
recent work of Donaldson, Wang, and  Phong-Sturm, we show that the
statement holds for higher rank vector bundles over compact
algebraic manifolds of arbitrary dimension that admit constant
scalar curvature metric and have discrete automorphism group.

\end{abstract}
\maketitle
\pagenumbering{arabic}

\newtheorem*{Mthm}{Main Theorem}
\newtheorem{Thm}{Theorem}[section]
\newtheorem{Prop}[Thm]{Proposition}
\newtheorem{Lem}[Thm]{Lemma}
\newtheorem{Cor}[Thm]{Corollary}
\newtheorem{Def}[Thm]{Definition}
\newtheorem{Guess}[Thm]{Conjecture}
\newtheorem{Ex}[Thm]{Example}
\newtheorem{Rmk}{Remark}
\newtheorem{Not}{Notation}
\def\thesection{\arabic{section}}
\renewcommand{\theThm} {\thesection.\arabic{Thm}}

\section{Introduction}

 A central notion in geometric invariant theory (GIT) is the
concept of stability. Stability plays a significant role in
forming quotient spaces of projective varieties for which
geometric invariant theory was invented. One can define
Mumford-Takemoto slope stability for holomorphic vector bundles,
and also there is a notion of Gieseker stability which is more in
the realm of geometric invariant theory. It is well-known that
over algebraic curves, all of these different notions coincide. It
was known from the work of Narasimhan and Seshadri that a
holomorphic vector bundle over a compact Riemann surface is
poly-stable if and only the bundle admits a projectively flat
connection. The picture became complete with the later work of
Donaldson, Uhlenbeck and Yau (\cite{D1}, \cite{D2} \cite{UY}).
They proved that over a compact K\"ahler manifold, a holomorphic
vector bundle is poly-stable if and only if it admits a
Hermitian-Einstein metric. This is known as the Hitchin-Kobayashi
correspondence. By a conjecture of Yau, one would also expect such
a correspondence for polarized algebraic manifolds. In other
words, the existence of extremal metrics on such a manifold should
be equivalent to being stable in some GIT sense. In \cite{Zh},
Zhang introduced the concept of balanced embedding and proved that
the existence of balanced embedding of a polarized algebraic
variety is equivalent to stability of Chow point of the variety.
Zhang's result has been reproven by Lu in \cite{L} and Phong and
Sturm in \cite{PS1}. The same correspondence was proven for vector
bundles by Wang in \cite{W1}. Later in \cite{D3}, Donaldson proved
that the existence of constant scalar curvature K\"ahler metrics
implies existence of balanced metrics and hence asymptotic Chow
stability. The converse is not yet known.

Earlier, in \cite{M}, Morrison proved that for the
projectivization of a rank two holomorphic vector bundle over a
compact Riemann surface, Chow stability is equivalent to the
stability of the bundle. Using ideas from the recent research
discussed above, in this article we generalize one direction of
Morrison's result for higher rank vector bundles over compact
algebraic manifolds of arbitrary dimension that admit constant
scalar curvature metric and have discrete automorphism group.

To state the precise result, let $X$ be a compact complex manifold
of dimension $m$ and $\pi : E\rightarrow X$ a holomorphic vector
bundle of rank $r$ with dual bundle $E^*$. This gives a
holomorphic fibre bundle $\mathbb{P}E^*$ over $X$ with fibre
$\mathbb{P}^{r-1}$. One can pull back the vector bundle $E$ to
$\mathbb{P}E^*$. We denote the tautological line bundle on
$\mathbb{P}E^*$ by $\mathcal{O}_{\mathbb{P}E^*}(-1)$ and its dual
by $\mathcal{O}_{\mathbb{P}E^*}(1)$. Let $L \rightarrow X$ be an
ample line bundle on $X$ and $\omega \in 2\pi  c_{1}(L)$ be a
K\"ahler form. Since $L$ is ample, there is an integer $k_{0} $ so
that for any $k\geq k_{0}$, $\mathcal{O}_{\mathbb{P}E_{k}^*}(1)$
is very ample over $\mathbb{P}E^*$, where $E_{k}=E\otimes
L^{\otimes k}$. Note that $\mathbb{P}E_{k}^* \cong \mathbb{P}E^*$
and $\mathcal{O}_{\mathbb{P}E_{k}^*}(1)\cong
\mathcal{O}_{\mathbb{P}E^*}(1) \otimes \pi^* L^{k}$. The theorem
we shall prove is the following:

\begin{thm}\label{thm1}
Suppose that $Aut(X)$ is discrete and $X$ admits a constant scalar
curvature K\"ahler metric in the class of $2\pi  c_{1}(L)$. If $E$
is Mumford stable, then there exists $k_{0}$ such that
$$(\mathbb{P}E^*,\mathcal{O}_{\mathbb{P}E^*}(1)\otimes \pi^* L^k)$$
is Chow stable for $k \geq k_{0}$.

\end{thm}

One of the earliest results in this spirit is the work of Burns
and De Bartolomeis in \cite{BD}. They construct a ruled surface
which does not admit any extremal metric in certain cohomology
class. In \cite{H1}, Hong proved that there are constant scalar
curvature K\"ahler metrics on the projectivization of stable
bundles over curves. In \cite{H2} and \cite{H3}, he generalizes
this result to higher dimensions with some extra assumptions.
Combining Hong's results with Donaldson's,
$(\mathbb{P}E^*,\mathcal{O}_{\mathbb{P}E_{m}^*}(n))$ is Chow
stable for $m,n \gg 0$ when the bundle $E$ is stable. Note that it
differs from our result, since it implies the Chow stability of
$(\mathbb{P}E^*,\mathcal{O}_{\mathbb{P}E_{m}^*}(n))$ for $n$ big
enough.

In \cite{RT}, Ross and Thomas developed the notion of slope
stability for polarized algebraic manifolds. As one of the
applications of their theory, they proved that if
$(\mathbb{P}E^*,\mathcal{O}_{\mathbb{P}E^*}(1)\otimes \pi^* L^k)$
is slope semi-stable for $k \gg 0$, then $E$ is a slope semistable
bundle and $(X,L)$ is a slope semistable manifold.Again note that
they look at stability of $\mathbb{P}E^*$ with respect to
polarizations $\mathcal{O}_{\mathbb{P}E_{m}^*}(n)$ for $n$ big
enough. For the case of one dimensional base, however they showed
stronger results. In this case they proved that if
$(\mathbb{P}E^*,\mathcal{O}_{\mathbb{P}E^*}(1)\otimes \pi^* L)$ is
slope (semi, poly) stable for any ample line bundle $L$, then $E$
is a slope (semi, poly) stable bundle.

In order to prove Theorem \ref{thm1} we use the concept of
\emph{balanced metrics} (See Definition \ref{ndef1}). Combining
the results of Luo, Phong, Sturm and Zhang on the relation between
balanced metrics and stability, it suffices to prove the following

\begin{thm}\label{thm2}

Let $X$ be a compact complex manifold and $L \rightarrow X$ be an
ample line bundle. Suppose that $X$ admits a constant scalar
curvature K\"ahler metric in the class of $2\pi  c_{1}(L)$ and
$Aut(X)$ is discrete. Let $E\rightarrow X$ be a holomorphic vector
bundle on $X$. If $E $ is Mumford stable, then
$\mathcal{O}_{\mathbb{P}E^*}(1)\otimes \pi^* L^k$ admits balanced
metrics for $k \gg 0$.

\end{thm}
The balanced condition may be formulated in terms of Bergman
kernels. First, we show that there exists an asymptotic expansion
for the Bergman kernel of
$(\mathbb{P}E^*,\mathcal{O}_{\mathbb{P}E^*}(1)\otimes \pi^* L^k)$.
Fix a positive hermitian metric $\sigma$ on $L$ such that $Ric
(\sigma)=\omega$. For any hermitian metric $g$ on
$\mathcal{O}_{\mathbb{P}E^*}(1)$, we define the sequence of volume
forms $d\mu_{g,k}$ on $\mathbb{P}E^*$ as follows
$$d\mu_{g,k}=k^{-m}\frac{(\omega_{g} + k
\pi^*\omega)^{m+r-1}}{(m+r-1)!}= \sum_{j=0}^{m}
 k^{j-m} \frac{\omega_{g}^{m+r-1-j}}{(m+r-j)!} \wedge
 \frac{\pi^*\omega^j}{j!},$$ where $\omega_{g}=Ric(g)$.

Let $\rho_{k}(g,\omega)$ be the Bergman kernel of
$H^{0}(\mathbb{P}E^*,\mathcal{O}_{\mathbb{P}E^*}(1)\otimes
\pi^*L^{k})$ with respect to the $L^2$-inner product $L^2(g
\otimes \sigma^{\otimes k}, d\mu_{k,g})$. We prove the
 following

\begin{thm}{\label{thmH1}}

For any hermitian metric $h$ on $E$ and K\"ahler form $\omega \in
2\pi c_{1}(L)$, there exist smooth endomorphisms
$\widetilde{B}_{k}(h,\omega)$ such that $$\rho_{k}(g,\omega)([v])=
C_{r}^{-1}tr \big( \lambda(v,h)\widetilde{B}_{k}(h,\omega)
\big),$$ where $g$ is the Fubini-Study metric on
$\mathcal{O}_{\mathbb{P}E^*}(1)$ induced by the hermitian metric
$h$. Moreover,
\begin{enumerate}
\item There exist smooth endomorphisms  $A_{i}(h,\omega) \in
\Gamma(X,E)$ such that the following asymptotic expansion holds as
$k \longrightarrow \infty$,
$$ \widetilde{B}_{k}(h,\omega) \sim
k^m+A_{1}(h,\omega)k^{m-1}+\dots.$$

\item In particular
$$A_{1}(h,\omega)= \frac{i}{2\pi} \Lambda F_{(E,h)}-
\frac{i}{2\pi r}tr( \Lambda F_{(E,h)})I_{E}+ \frac{(r+1)}{2r}
S(\omega) I_{E},$$ where $\Lambda $ is the trace operator acting
on $(1,1)$-forms with respect to the K\"ahler form $\omega$ and
$F_{(E,h)}$ is the curvature of $(E,h)$ and $S(\omega)$ is the
scalar curvature of $\omega$.

\item The asymptotic expansion holds in $C^{\infty}$. More
precisely, for any positive integers $a$ and $p$, there exists a
positive constant $K_{a,p,\omega,h}$ such that
$$\Big |\big|\widetilde{B}_{k}(h,\omega)-\big(
k^m+\dots+A_{p}(h,\omega)k^{m-p} \big)\big|\Big|_{C^a}\leq
K_{a,p,\omega,h} k^{m-p-1}.$$ Moreover the expansion is uniform in
the sense that there exists a positive integer $s$ such that if
$h$ and $\omega$ run in a bounded family in $C^s$ topology and
$\omega$ is bounded from below, then the constants $K_{a,p,\omega,h}$
are bounded by a constant depending only on $a$ and $p$.

\end{enumerate}

\end{thm}

Finding balanced metrics on $\mathcal{O}_{\mathbb{P}E^*}(1)\otimes
\pi^* L^k$ is basically the same as finding solutions to the
equations $\rho_{k}(g,\omega)= \textrm{Constant}.$ Therefore in
order to prove Theorem \ref{thm2}, we need to solve the equations
$\rho_{k}(g,\omega)= \textrm{Constant}$ for $k \gg 0$. Now if
$\omega$ has constant scalar curvature and $h$ satisfies the
Hermitian-Einstein equation $\Lambda_{\omega}F_{(E,h)}=\mu I_{E}$,
then $A_{1}(h,\omega)$ is constant. Notice that in order to make
$A_{1}$ constant, existence of Hermitian-Einstein metric is not
enough. We need the existence of constant scalar curvature
K\"ahler metric as well. Next, the crucial fact is that the
linearization of $A_{1}$ at $(h,\omega)$ is surjective. This
enables us to construct formal solutions as power series in
$k^{-1}$ for the equation $\rho_{k}(g,\omega)= \textrm{Constant}.$
Therefore, for any positive integer $q$, we can construct a
sequence of metrics $g_{k}$ on
$\mathcal{O}_{\mathbb{P}E^*}(1)\otimes \pi^*L^{k}$ and bases
$s^{(k)}_{1},...,s_{N}^{(k)}$ for
$H^0(\mathbb{P}E^*,\mathcal{O}_{\mathbb{P}E^*}(1))$ such that
$$\sum |s^{(k)}_{i}|_{g_{k}}^2=1,$$
$$\int \langle s^{(k)}_{i},s^{(k)}_{j}\rangle_{g_{k}} dvol_{g_{k}}=
D_{k}I+M_{k},$$ where $D_{k} \rightarrow C_{r}$ as $k \rightarrow
\infty$ (See \eqref{eq16} for definition of $C_{r}$.), and $M_{k}$
is a trace-free hermitian matrix such that
$||M_{k}||_{\textrm{op}}=o(k^{-q-1})$ as $k\rightarrow \infty$.

The next step is to perturb these almost balanced metrics to get
balanced metrics. As pointed out by Donaldson, the problem of
finding balanced metric can be viewed also as a finite dimensional
moment map problem solving the equation $M_{k}=0$. Indeed,
Donaldson shows that $M_{k}$ is the value of a moment map
$\mu_{D}$ on the space of ordered bases with the obvious action of
$SU(N)$. Now, the problem is to show that if for some ordered
basis $\underline{s}$, the value of moment map is very small, then
we can find a basis at which moment map is zero. The standard
technique is flowing down $\underline{s}$ under the gradient flow
of $|\mu_{D}|^2$ to reach a zero of $\mu_{D}$. We need a
Lojasiewicz type inequality to guarantee that the flow converges
to a zero of the moment map. We do this in Section 3 by adapting
Phong-Sturm proof to our situation.

Here is the outline of the paper: In Section $2$, we review
Donaldson's moment map setup. We follow Phong and Sturm treatment
from  (\cite {PS2}). In Section $3$, we obtain a lower bound for
the derivative of the moment map by adapting the argument in
(\cite{PS2}) to our setting. In Section $4$, we show how to
perturb almost balanced metrics to obtain balanced metrics in the
general setting of Section $3$. In order to do that, we use the
estimate obtained in Theorem \ref{thm4} to apply the Donaldson's
version of inverse function theorem(Proposition \ref{prop1}). In
Section $5$, we prove the existence of an asymptotic expansion for
the Bergman kernel of $\mathcal{O}_{\mathbb{P}E^*}(1)\otimes
\pi^*L^{k}$ using results of Catlin and Zelditch. Section $6$ is
devoted to constructing almost balanced metrics on
$\mathcal{O}_{\mathbb{P}E^*}(1)\otimes \pi^*L^{k}$  using the
asymptotic expansion obtained in Section $5$.

\thanks{\textbf{Acknowledgements:} I am sincerely grateful to Richard Wentworth for introducing me the subject and many helpful discussions and suggestions
on the subject and his continuous help, support and encouragement.
I would also like to thank Bo Berndtsson, Hamid Hezari, Duong Hong
Phong , Julius Ross and Steve Zelditch for many helpful
discussions and suggestions.}

\section{Moment Map Setup}

In this section, we review Donaldson's moment map setup. We follow
the notation of \cite{PS2}.

Let $(Y,\omega_{0})$ be a compact K\"ahler manifold of dimension
$n$ and $\mathcal{O}(1)\rightarrow Y$ be a very ample line bundle
on $Y$ equipped with a Hermitian metric $g_{0}$ such that
$Ric(g_{0})=\omega_{0}$. Since $\mathcal{O}(1)$ is very ample,
using global sections of $\mathcal{O}(1)$, we can embed $Y$ into
$\mathbb{P}(H^0(Y,\mathcal{O}(1))^*)$. A choice of ordered basis
$\underline{s}=(s_{1},...,s_{N})$ of $H^0(Y,\mathcal{O}(1))$ gives
an isomorphism between $\mathbb{P}(H^0(Y,\mathcal{O}(1))^*)$ and
$\mathbb{P}^{N-1}$. Hence for any such $\underline{s}$, we have an
embedding $\iota_{\underline{s}}:Y\hookrightarrow
\mathbb{P}^{N-1}$ such that $\iota_{\underline{s}}^*
\mathcal{O}_{\mathbb{P}^N}(1)=\mathcal{O}(1)$. Using
$\iota_{\underline{s}}$, we can pull back the Fubini-Study metric
and K\"ahler form of the projective space to $\mathcal{O}(1)$ and
$Y$ respectively.
\begin{Def}\label{ndef1}
An embedding $\iota_{\underline{s}}$ is called balanced if
$$\int_{Y} \langle  s_{i}, s_{j}
\rangle_{\iota_{\underline{s}}^*h_{\textrm{FS}}}\frac{\iota_{\underline{s}}^*\omega_{\textrm{FS}}}{n!}=\frac{V}{N}\delta_{ij},$$
where $V=\int_{Y}\frac{\omega_{0}^n}{n!}$. A hermitian
metric(respectively a K\"ahler form) is called balanced if it is
the pull back $\iota^*_{\underline{s}}h_{\textrm{FS}}$
(respectively $\iota^*_{\underline{s}}\omega_{\textrm{FS}}$) where
$\iota_{\underline{s}}$ is a balanced embedding.
\end{Def}

There is an action of $SL(N)$ on the space of ordered bases of
$H^0(Y,\mathcal{O}(1)).$ Donaldson defines a symplectic form on
the space of ordered bases of $H^0(Y,\mathcal{O}(1))$  which is
invariant under the action of $SU(N)$. So there exists an
equivariant moment map on this space such that its zeros are
exactly balanced bases.

More precisely we define
$$\widetilde{\mathcal{Z}}=\{\underline{s}=(s_{1},...,s_{N})|s_{1},...,s_{N}
\textrm{ a basis of } H^0(Y,\mathcal{O}(1))\}/\mathbb{C}^*$$ and
$\mathcal{Z}=\widetilde{\mathcal{Z}}/\mathbb{P}Aut(Y,\mathcal{O}(1))$.
Donaldson defines a symplectic form  $\Omega_{D}$ on
$\mathcal{Z}$. There is a natural action of $SU(N)$ on
$(\mathcal{Z},\Omega_{D})$ which preserves the symplectic form
$\Omega_{D}$. The moment map for this action is defined by
$$\mu_{D}(\underline{s})=i[ \langle s_{\alpha}, s_{\beta}
\rangle_{h_{\underline{s}}}-\frac{V}{N} \delta_{\alpha,\beta}],
$$ where $h_{\underline{s}}$ is the $L^2$- inner product with respect to the pull back of Fubini-Study metric and Fubini-Study K\"ahler form via the embedding $\iota_{\underline{s}}$.
Also we identify  $su(N)^*$ with $su(N)$ using the invariant inner
product on $su(N)$, where $su(N)$ is the Lie algebra of the group
$SU(N)$ and $su(N)^*$ is its dual. (For construction of
$\Omega_{D}$ and more details see (\cite{D3}) and (\cite{PS2}) .)

Using Deligne's pairing, Phong and Sturm construct another
symplectic form on $\mathcal{Z}$ as follows:

Let $$\widetilde{\mathcal{Y}}=\{(x,\underline{s}) | x \in
\mathbb{P}^{N-1}, \underline{s}=(s_{1},...,s_{N}), x \in
\iota_{\underline{s}}(Y)\}$$ and
$\mathcal{Y}=\widetilde{\mathcal{Y}}/\mathbb{P}Aut(Y,\mathcal{O}(1))$.
One obtains a holomorphic fibration $\mathcal{Y}\rightarrow
\mathcal{Z}$ where every fibre is isomorphic to $Y$. Let
$p:\mathcal{Y}\rightarrow \mathbb{P}^{N-1} $ be the projection on
the first factor. Then define a hermitian line bundle
$\mathcal{M}$ on $\mathcal{Z}$ by $$\mathcal{M}=\langle
p^*\mathcal{O}_{\mathbb{P}^{N-1}}(1),...,p^*\mathcal{O}_{\mathbb{P}^{N-1}}(1)
\rangle (\frac{\mathcal{Y}}{\mathcal{Z}})$$ which is the Deligne's
pairing of $(n+1)$ copies of
$p^*\mathcal{O}_{\mathbb{P}^{N-1}}(1)$. Denote the curvature of
this hermitian line bundle by $\Omega_{\mathcal{M}}$. It follows
from properties of Deligne's pairing that

\begin{equation}\label{eq1}\Omega_{\mathcal{M}}= \int_{\mathcal{Y}/\mathcal{Z}}
\omega_{\textrm{FS}}^{n+1}.\end{equation}

Since $SU(N)$ is semisimple, there is a unique equivariant moment
map $\mu_{\mathcal{M}}: \mathcal{Z}\rightarrow su(N)$ for the
action of $SU(N)$ on $(\mathcal{Z},\Omega_{\mathcal{M}})$.

\begin{thm}\label{thm3}(\cite[Theorem 1]{PS2})
$\Omega_{\mathcal{M}}=\Omega_{D}$ and $\mu_{\mathcal{M}}=\mu_{D}$.

\end{thm}
Let $\xi$ be an element of the Lie algebra $su(N)$. Since $SU(N)$
acts on $\mathcal{Z}$, the infinitesimal action of $\xi$ defines a
vector field $\sigma_{\mathcal{Z}}(\xi)$ on $\mathcal{Z}$. Fixing
a point $z \in \mathcal{Z}$, we have a linear map $\sigma_{z}
:su(N)\rightarrow T_{z}\mathcal{Z}$. Let $\sigma_{z}^* $ be its
adjoint with respect to the metric on $T\mathcal{Z}$ and the
invariant metric on $su(N)$. Then we get the operator
$$Q_{z}=\sigma_{z}^* \sigma_{z}: su(N)\rightarrow su(N).$$
Define $\Lambda_{z}^{-1}$ as the smallest eigenvalue of $Q_{z}$.
In \cite{D3}, Donaldson proves the following.

\begin{prop}\label{prop1}(\cite[Proposition 17]{D3})
Suppose given $z_{0} \in \mathcal{Z}$ and real numbers $\lambda,
\delta$ such that for all $z=e^{i\xi}z_{0} $ with
$|\xi|\leq\delta$ and $\xi \in su(N)$, $\Lambda_{z}\leq\lambda$.
Suppose that $\lambda |\mu(z_{0})|\leq \delta$, then there exists
$w=e^{i\eta}$ with $\mu(w)=0$, where $|\eta|\leq\lambda
|\mu(z_{0})|.$

\end{prop}

\section{Eigenvalue Estimates}

In this section, we obtain a lower bound for the derivative of the
moment map $\mu_{D}$. This is equivalent to an upper bound for the
quantity $\Lambda_{z}$ introduced in the previous section. In
order to do this, we adapt the argument of Phong and Sturm  to our
setting. The main result is Theorem \ref{thm4}.

Let $(Y,\omega_{0})$  and $\mathcal{O}(1)\rightarrow Y$ be as in
the previous section. Let $(L, h_{\infty})$ be a Hermitian line
bundle over $Y$ such that $\omega_{\infty}=Ric(h_{\infty})$ is a
semi positive $(1,1)$-form on $Y$. Define
$\widetilde{\omega}_{0}=\omega_{0}+k\omega_{\infty}$. For the rest
of this section and next section let $m$ be the smallest integer
such that $\omega_{\infty}^{m+1}=0.$ Also assume that
$\omega_{0}^{n-m}\wedge \omega_{\infty}^m$ is a volume form and
there exist positive constant $n_{1}$ and $n_{2}$ such that

\begin{align}\label{eq4}  &N_{k}= \dim H^0(Y,\mathcal{O}(1)\otimes
L^{k})=n_{1}k^m+O(k^{m-1}).\\ \label{eq14}
&V_{k}=\int_{Y}(\omega_{0}+k\omega_{\infty})^n=n_{2}k^{m}+O(k^{m-1}).
\end{align}
Notice that \eqref{eq14} is implied from the fact that
$\omega_{0}^{n-m}\wedge \omega_{\infty}^m$ is a volume form and
$\omega_{\infty}^{m+1}=0$.

The case important for this paper is the following:

\begin{Ex}

Let $(X,\omega_{\infty})$ be a compact K\"ahler manifold of
dimension $m$ and $L$ be a very ample holomorphic line bundle on
$X$ such that $\omega_{\infty} \in 2\pi c_{1}(L)$. Let $E$ be a
holomorphic vector bundle on $X$ of rank $r$ such that the line
bundle $\mathcal{O}_{\mathbb{P}E^*}(1)\rightarrow Y=\mathbb{P}E^*$
is an ample line bundle. We denote the pull back of
$\omega_{\infty}$ to $\mathbb{P}E^*$ by $\omega_{\infty}$. Then
$\omega_{\infty}^{m+1}=0$ and by Riemann-Roch formula we have
$$\dim H^0(Y,\mathcal{O}(1)\otimes L^{k})=\dim H^0(X,E\otimes
L^{k})=\frac{r}{m!}\int_{X}c_{1}(L)^m k^m+O(k^{m-1}).$$

\end{Ex}

The following lemma is clear.
\begin{lem}\label{lem5}
Let $h_{k} $ be a sequence of hermitian metrics  on
$\mathcal{O}(1)\otimes L^{k}$ and let
$\underline{s}^{(k)}=(s_{1}^{(k)},...,s_{N}^{(k)})$ be a sequence
of ordered bases for $H^0(Y,\mathcal{O}(1)\otimes L^{k})$. Suppose
that for any $k$  $$\sum |s_{i}^{(k)}|_{h_{k}}^2=1$$ and
$$\int_{Y}\langle s_{i}^{(k)}, s_{j}^{(k)} \rangle
_{h_{k}}dvol_{h_{k}}=D^{(k)}\delta_{ij}+M^{(k)}_{ij},$$
 where $D^{(k)}$ is a scalar  and $M^{(k)}$ is a trace-free hermitian
 matrix. Then $$D^{(k)}=\frac{V_{k}}{N_{k}} \rightarrow \frac{n_{2}}{n_{1}}  \,\,\,\,\,\textrm{ as} \,\,\, k\rightarrow \infty, $$
where the constants $n_{1}$ and $n_{2}$ are defined by \eqref{eq4}
and \eqref{eq14}.

\end{lem}

We start with the notion of $R$-boundedness introduced originally
by Donaldson in \cite{D3}.

\begin{Def}\label{def1}

Let $R$ be a real number with $R >1$ and $a\geq 4$ be a fixed
integer and let $\underline{s}=(s_{1},...,s_{N})$ be an ordered
basis for $H^0(Y,\mathcal{O}(1)\otimes
 L^{k})$. We say $\underline{s}$ has $R$-bounded geometry if the
 K\"ahler form
  $\widetilde{\omega}=\iota^*_{\underline{s}}\omega_{\textrm{FS}}$ satisfies the following conditions
 \begin{itemize}

\item

$||\widetilde{\omega}-\widetilde{\omega}_{0}||_{C^{a}(\widetilde{\omega}_{0})}\leq
R$, where $\widetilde{\omega}_{0}=\omega_{0}+ k \omega_{\infty}$.

\item$\widetilde{\omega} \geq \frac{1}{R} \widetilde{\omega}_{0}.$

 \end{itemize}

\end{Def}

 Recall the definition of $\Lambda_{z}$ from the previous section. The main result of this section is the following.

 \begin{thm}\label{thm4}

 Assume $Y$ does not have any nonzero holomorphic vector fields. For any $R>1$, there are positive constants $C$ and  $\epsilon \leq n_{2}/10n_{1}$
 such that, for any $k$, if the basis $\underline{s}=(s_{1},...,s_{N})$ of $H^0(Y,\mathcal{O}(1)\otimes
 L^{k})$ has $R$-bounded geometry, and if $||\mu_{D}(\underline{s})||_{\textrm{op}}\leq\epsilon
 $, then $$\Lambda_{\underline{s}} \leq Ck^{2m+2}.$$

\end{thm}

The rest of this section is devoted to the proof of Theorem
\ref{thm4}. Notice that the estimate $\Lambda_{z}\leq Ck^{2m+2}$
is equivalent to the estimate
\begin{equation}\label{eq5}|\sigma_{\mathcal{Z}}(\xi)|^2\geq ck^{-(2m+2)}||\xi||^2.\end{equation}

On the other hand \eqref{eq1} and Theorem \ref{thm3} imply that
\begin{equation}\label{eq6}|\sigma_{\mathcal{Z}}(\xi)|^2=  \int_{Y}\iota_{Y_{\xi},\overline{Y_{\xi}}} \omega_{\textrm{FS}}^{n+1}.\end{equation}
Hence, in order to establish Theorem \ref{thm4}, we need to
estimate the quantity $ \int_{Y} \iota_{Y_{\xi},
\overline{Y_{\xi}}} \omega_{\textrm{FS}}^{n+1}$ from below.

For the rest of this section, fix an ordered basis
$\underline{s}^{(k)}=(s_{1},...,s_{N})$ of
$H^0(Y,\mathcal{O}(1)\otimes L^{k})$ and let
$M^{(k)}=-i\mu_{D}(\underline{s}^{(k)})$. It gives an embedding
$\iota=\iota_{\underline{s}^{(k)}}:Y\rightarrow \mathbb{P}^{N-1}$,
where $N=N_{k}=\dim H^0(Y,\mathcal{O}(1)\otimes L^{k})$. For any
$\xi \in su(N)$, we have a vector field $Y_{\xi}$ on
$\mathbb{P}^{N-1}$ generated by the infinitesimal action of $\xi$.

Every tangent vector to $\mathbb{P}^{N-1}$ is given by pairs
$(z,v)$ modulo an equivalence relation $\sim$ . This relation is
defined as follows:
$$(z,v)\sim(z',v') \textrm{ if } z'=\lambda z \textrm{ and }
v'-\lambda v=\mu z \textrm{ for some } \lambda \in \mathcal{C}^*
\textrm{ and } \mu \in \mathbb{C}.$$ For a tangent vector
$[(z,v)]$, the Fubini-Study metric is given by
$$||[(z,v)]||^2=\frac{v^*vz^*z-(z^*v)^2}{(z^*z)^2}.$$ Since the
vector field $Y_{\xi}$ is given by $[{z,\xi z}]$, we have

\begin{equation}\label{eq7}||Y_{\xi}(z)||^2= \frac{-(z^*\xi z)^2+
(z^* \xi^2 z)(z^* z)}{(z^* z)^2}.\end{equation} We have the
following exact sequence of vector bundles over $Y$
$$0\rightarrow TY \rightarrow \iota^* T\mathbb{P}^{N-1}
\rightarrow Q \rightarrow 0.$$

Let $\mathcal{N}\subset \iota^* T\mathbb{P}^{N-1} $ be the
orthogonal complement of $TY$. Then as smooth vector bundles, we
have
$$\iota ^* T\mathbb{P}^{N-1}= TY\oplus \mathcal{N}.$$ We denote the
projections onto the first and second component by $\pi_{T}$ and
$\pi_{\mathcal{N}}$ respectively. Define
 $$\sigma_{t}(z)=\exp (it \xi) z,$$
 $$\varphi_{t}(z)=\log \frac{|\sigma_{t}(z)|}{|z|}.$$

Direct computation shows that

 \begin{equation}\label{eq8}\frac{d}{dt}\Big|_{t=0}
\varphi_{t}(z)=2i\frac{z^*\xi z}{z^*z},\end{equation}

 \begin{equation}\label{eq9}\frac{d^2}{dt^2}\Big|_{t=0} \varphi_{t}(z)=4\frac{(z^*\xi z)^2-
(z^* \xi^2 z)(z^* z)}{(z^* z)^2}.\end{equation} The following is
straightforward.

\begin{prop}\label{prop2}

For any $\xi \in su(N)$, we have
$$||\pi_{\mathcal{N}}Y_{\xi}||^2_{L^2(Y,TY)}= \int_{Y}
\iota_{Y_{\xi}, \overline{Y_{\xi}}} \omega_{\textrm{FS}}^{n+1}$$

\end{prop}

Therefore, the estimate in Theorem \ref{thm4} will follow from:

\begin{equation}\label{eq10}||\xi||^2\leq c_{R} k^m||Y_{\xi}||^2 \end{equation}
\begin{equation}\label{eq11} c^{'}_{R}||\pi_{T}Y_{\xi}||^2\leq k^{m+2}||\pi_{\mathcal{N}}Y_{\xi}||^2\end{equation}
\begin{equation}\label{eq12}||Y_{\xi}||^2 =||\pi_{T}Y_{\xi}||^2+||\pi_{\mathcal{N}}Y_{\xi}||^2\end{equation}

We will prove \eqref{eq10} in Proposition \ref{prop3} and
\eqref{eq11} in Proposition \ref{prop5}. Assuming these, we give
the Proof of Theorem \ref{thm4}.

\begin{proof}[Proof of Theorem \ref{thm4}.]
By \eqref{eq6}, we have $$|\sigma_{\mathcal{Z}}(\xi)|^2=
\int_{Y}\iota_{Y_{\xi},\overline{Y_{\xi}}}
\omega_{\textrm{FS}}^{n+1}.$$ Applying Proposition \ref{prop2}, we
get
$$|\sigma_{\mathcal{Z}}(\xi)|^2=||\pi_{\mathcal{N}}Y_{\xi}||^2.$$
Thus, in order to prove Theorem \ref{thm4}, we need to show that
$$||\pi_{\mathcal{N}}Y_{\xi}||^2 \geq c_{R}k^{-(m+3)}||\xi||^2.$$
By \eqref{eq10}, we have $$||\xi||^2\leq c_{R} k^m||Y_{\xi}||^2=
c_{R}
k^m||\pi_{\mathcal{N}}Y_{\xi}||^2+c_{R}k^m||\pi_{T}Y_{\xi}||^2.$$
Hence \eqref{eq11} implies that \begin{align*}||\xi||^2&\leq c_{R}
k^m||\pi_{\mathcal{N}}Y_{\xi}||^2+c_{R}c_{R}'
k^{2m+2}||\pi_{\mathcal{N}}Y_{\xi}||^2\\& \leq
c_{R}''k^{2m+2}||\pi_{\mathcal{N}}Y_{\xi}||^2.\end{align*}

\end{proof}

\begin{lem}\label{lem1}
There exists a positive constant $c$ independent of $k$ such that
for any  $f \in C^{\infty}(Y)$, we have
$$ c\int_{Y} f^2 \widetilde{\omega}_{0}^n \leq k^m \int_{Y}
\overline{\partial}f\wedge \partial f \wedge
\widetilde{\omega}_{0}^{n-1} +k^{-m}\Big( \int_{Y} f
\widetilde{\omega}_{0}^n \Big)^2$$

\end{lem}

\begin{proof}

In the proof of this Lemma, we put
$\omega_{k}=\omega_{0}+k\omega_{\infty}$ and
$\alpha=\omega_{1}=\omega_{0}+\omega_{\infty}.$ For $k \geq 1,$ we
have
$$k^{-m}\omega_{k}^n \leq  \alpha^n \leq \omega_{k}^n. $$

Assume that the statement is false. So, there exists a subsequence
$k_{j} \rightarrow \infty$  and a sequence of functions $f_{j}$
such that $\int_{Y} f_{j}^2 \omega_{k_{j}}^n=1 $ and
$$k^m \int_{Y}\overline{\partial}f_{j}\wedge \partial f_{j} \wedge
\omega_{k_{j}}^{n-1} +k_{j}^{-m}\Big( \int_{Y}
f_{j}\omega_{k_{j}}^n \Big)^2 \rightarrow 0  $$ as $k\rightarrow
\infty.$ We define $||f||^2=\int_{Y} f^2 \alpha^n$. Hence
$$||f_{j}||^2=\int_{Y} f_{j}^2 \alpha^n \geq k_{j}^{-m}\int_{Y}
f_{j}^2 \omega_{k_{j}}^n =k_{j}^{-m}.$$

Let $g_{j}=f_{j} \big/||f_{j}||.$ We have
\begin{align*}\int_{Y}|\partial g_{j}|_{\alpha}^2 \alpha^n &=
\int_{Y}\overline{\partial}g_{j}\wedge \partial g_{j} \wedge
\alpha^{n-1}\\&=||f_{j}||^{-2}\int_{Y}\overline{\partial}f_{j}\wedge\partial
f_{j}\wedge\alpha^{n-1}\\&\leq
k_{j}^m\int_{Y}\overline{\partial}f_{j}\wedge\partial
f_{j}\wedge\omega_{k_{j}}^{n-1}\rightarrow 0 \,\,\,\,\,\, \textrm{
as }\,\,\,\, k\rightarrow \infty .\end{align*} On the other hand
$\int_{Y} g_{j}^2 \alpha^n=1$ which implies that the sequence
${g_{j}}$ is bounded in $L^2_{1}(\alpha^n)$. Hence, ${g_{j}}$ has
a subsequence which converges in $L^2(\alpha^n)$ and converges
weakly in $L^2_{1}(\alpha^n)$ to a function $g \in
L^2_{1}(\alpha^n)$. Without loss of generality, we can assume that
the whole sequence converges. Since $\int_{Y}|\partial
g_{j}|_{\alpha}^2 \alpha^n \rightarrow 0$ as $k\rightarrow
\infty$, it can be easily seen that $g$ is a constant function. We
have \begin{align*}k_{j}^{-m}\big|  \int_{Y} (g_{j}-g)
\omega_{k_{j}}^n \big|&\leq k_{j}^{-m}\int_{Y} |g_{j}-g|
\omega_{k_{j}}^n \\& \leq \int_{Y} |g_{j}-g| \alpha^n\\
&\leq C(\int_{Y} | g_{j}-g|^2 \alpha^n)^{\frac{1}{2}} \rightarrow
0,\end{align*} where $C^2= \int_{Y} \alpha^n$ does not depend on
$k$. Hence
$$k_{j}^{-m}\big|  \int_{Y} (g_{j}-g) \omega_{k_{j}}^n  \big| \rightarrow 0.$$
Since $g$ is a constant function and $
\int_{Y}\omega_{k_{j}}^n=n_{2}k_{j}^m+O(k_{j}^{m-1})$, we get
$$k_{j}^{-m}\int_{Y} g_{j} \omega_{k_{j}}^n \rightarrow n_{2}g,$$
where $n_{2}$ is defined by \eqref{eq14}. On the other hand
\begin{align*} \Big(k_{j}^{-m} \int_{Y} g_{j}\omega_{k_{j}}^n
\Big)^2&=k_{j}^{-2m}||f_{j}||^{-2}\Big( \int_{Y}
f_{j}\omega_{k_{j}}^n \Big)^2 \\&\leq k_{j}^{-m}( \int_{Y}
f_{j}\omega_{k_{j}}^n \Big)^2 \rightarrow 0\end{align*} which
implies $g\equiv 0$. It is a contradiction since $||g_{j}||=1$ and
$g_{j} \rightarrow g$ in $L^2(\alpha^n)$.

\end{proof}

The proof of the following lemma can be found in (\cite[p.
704]{PS2}). For the sake of completeness, we give the details.

\begin{lem}\label{lem2}
  There exists a positive constant $c_{R}$ independent of $k$ such that for any K\"ahler form $\widetilde{\omega}\in c_{1}(\mathcal{O}(1)\otimes
L^{k})$  having $R$-bounded geometry and any $f \in
C^{\infty}(Y)$, we have
 $$c_{R}\int_{Y} f^2 \widetilde{\omega}^n \leq k^m\int_{Y} \overline{\partial}f\wedge \partial f \wedge \widetilde{\omega}^{n-1} +k^{-m}\Big( \int_{Y} f \widetilde{\omega}^n
 \Big)^2.$$

\end{lem}

\begin{proof}

Since $\widetilde{\omega }$ has $R$-bounded geometry, we have
$$R^{-1}\widetilde{\omega}_{0} \leq \widetilde{\omega } \leq2R \widetilde{\omega}_{0}.$$
Therefore,
$$c(2R)^{-n} \int_{Y} f^2 \widetilde{\omega}^n \leq c\int_{Y} f^2 \widetilde{\omega}_{0}^n \leq k^m\int_{Y} \overline{\partial}f\wedge \partial f \wedge \widetilde{\omega}_{0}^{n-1} +k^{-m}\Big( \int_{Y} f \widetilde{\omega}_{0}^n \Big)^2.$$
On the other hand, there exists a unique function $\phi$ such that
$\widetilde{\omega}-\widetilde{\omega}_{0}=\partial\overline{\partial}\phi$
and $\int_{Y} \phi \widetilde{\omega}_{0}^n=0$. Hence,
$$\widetilde{\omega}^n-\widetilde{\omega}_{0}^n=\partial\overline{\partial}\phi \wedge \sum_{j=0}^{n-1}\widetilde{\omega}^j \wedge \widetilde{\omega}_{0}^{n-j-1}$$
We have, \begin{align*}\Big| \int f
(\widetilde{\omega}^n-\widetilde{\omega}_{0}^n)    \Big|&= \Big|
\int f  \partial\overline{\partial}\phi \wedge
\sum_{j=0}^{n-1}\widetilde{\omega}^j \wedge
\widetilde{\omega}_{0}^{n-j-1}   \Big|\\&=\Big| \int
\overline{\partial}f \wedge  \partial\phi \wedge
\sum_{j=0}^{n-1}\widetilde{\omega}^j \wedge
\widetilde{\omega}_{0}^{n-j-1}   \Big|\\&\leq \sum_{j=0}^{n-1}
\int |\overline{\partial}f|_{\widetilde{\omega}_{0}}|\partial
\phi|_{\widetilde{\omega}_{0}}\big(\frac{\widetilde{\omega}}{\widetilde{\omega}_{0}}\big)^p
\widetilde{\omega}_{0}^n \\&\leq n(2R)^n \int
|\overline{\partial}f|_{\widetilde{\omega}_{0}}|\partial
\phi|_{\widetilde{\omega}_{0}} \widetilde{\omega}_{0}^n\\&\leq
C_{1 } \Big (\int
|\overline{\partial}f|^2_{\widetilde{\omega}_{0}}
\widetilde{\omega}_{0}^n \Big)^{\frac{1}{2}} \Big (\int
|\overline{\partial}\phi|^2_{\widetilde{\omega}_{0}}
\widetilde{\omega}_{0}^n \Big)^{\frac{1}{2}} \\&= C_{1 }  \Big
(\int \overline{\partial}f \wedge \partial f \wedge
\widetilde{\omega}_{0}^{n-1} \Big)^{\frac{1}{2}} \Big (\int
|\overline{\partial}\phi|^2_{\widetilde{\omega}_{0}}
\widetilde{\omega}_{0}^n \Big)^{\frac{1}{2}} .\end{align*} We will
show that
$$\int |\overline{\partial}\phi|^2_{\widetilde{\omega}_{0}} \widetilde{\omega}_{0}^n \leq C_{2} k^{2m}.$$
Since
$\widetilde{\omega}-\widetilde{\omega}_{0}=\partial\overline{\partial}\phi$
and
$||\widetilde{\omega}-\widetilde{\omega}_{0}||_{C^a(\widetilde{\omega}_{0})}
\leq R$, we have
$||\partial\overline{\partial}\phi||_{C^a(\widetilde{\omega}_{0})}\leq
R$. This implies that
$$||\triangle_{\widetilde{\omega}_{0}}\phi||_{\infty}\leq R.$$
Applying Lemma \ref{lem1} to $\phi$, we get $$c\int_{Y} \phi^2
\widetilde{\omega}_{0}^n \leq k^m \int_{Y}
\overline{\partial}\phi\wedge
\partial \phi \wedge \widetilde{\omega}_{0}^{n-1}$$
On the other hand \begin{align*}\int_{Y}
|\overline{\partial}\phi|^2_{\widetilde{\omega}_{0}}\widetilde{\omega}_{0}^n&=\int_{Y}
\overline{\partial}\phi\wedge
\partial \phi \wedge \widetilde{\omega}_{0}^{n-1}= \Big|  \int_{Y} \phi \triangle_{\widetilde{\omega}_{0}}\phi  \widetilde{\omega}_{0}^n
\Big|\\&\leq \Big(  \int_{Y} \phi^2 \widetilde{\omega}_{0}^n
\Big)^{\frac{1}{2}} \Big(  \int_{Y}
|\triangle_{\widetilde{\omega}_{0}}\phi|^2
\widetilde{\omega}_{0}^n \Big)^{\frac{1}{2}}\\&\leq
c^{\frac{-1}{2}}k^{\frac{m}{2}}  \Big(  \int_{Y}
|\overline{\partial}\phi|^2_{\widetilde{\omega}_{0}}
\widetilde{\omega}_{0}^n \Big)^{\frac{1}{2}} \Big( R^2 \int_{Y}
\widetilde{\omega}_{0}^n \Big)^{\frac{1}{2}}\end{align*}

$$=C k^{m} \Big(  \int_{Y} |\overline{\partial}\phi|^2_{\widetilde{\omega}_{0}} \widetilde{\omega}_{0}^n
\Big)^{\frac{1}{2}}$$ Therefore,

$$\int |\overline{\partial}\phi|^2_{\widetilde{\omega}_{0}} \widetilde{\omega}_{0}^n \leq C_{2} k^{2m}.$$
So, we get $$\Big| \int
f(\widetilde{\omega}^n-\widetilde{\omega}_{0}^n)    \Big| \leq C
k^{m}\Big ( \int_{Y} \overline{\partial}f\wedge
\partial f \wedge \widetilde{\omega}_{0}^{n-1} \Big)^{\frac{1}{2}}$$
On the other hand $$\frac{1}{2}\Big( \int_{Y} f
\widetilde{\omega}_{0}^n \Big)^2 \leq \Big( \int_{Y} f
\widetilde{\omega}^n \Big)^2+ \Big( \int_{Y} f
\big(\widetilde{\omega}^n-\widetilde{\omega}_{0}^n \big)\Big)^2$$
Hence, \begin{align*}\widetilde{C}\int_{Y} f^2
\widetilde{\omega}^n &\leq k^m\int_{Y} \overline{\partial}f\wedge
\partial f \wedge \widetilde{\omega}_{0}^{n-1} +2k^{-m}\Big(\Big(
\int_{Y} f \widetilde{\omega}^n \Big)^2+ \Big( \int_{Y} f
\big(\widetilde{\omega}^n-\widetilde{\omega}_{0}^n \big)\Big)^2
\Big)\\&\leq k^m\int_{Y} \overline{\partial}f\wedge
\partial f \wedge \widetilde{\omega}_{0}^{n-1} +2k^{-m}\Big(
\int_{Y} f \widetilde{\omega}^n \Big)^2+C_{3}k^m \int_{Y}
\overline{\partial}f\wedge
\partial f \wedge \widetilde{\omega}_{0}^{n-1}\\&\leq C_{4}\Big(  k^m\int_{Y} \overline{\partial}f\wedge \partial f \wedge \widetilde{\omega}^{n-1} +k^{-m}\Big( \int_{Y} f \widetilde{\omega}^n \Big)^2  \Big).\end{align*}

\end{proof}

\begin{prop}\label{prop3}

There exists a positive constant $c_{R}$ such that for any $\xi
\in su(N)$, we have
$$ ||\xi||^2\leq c_{R} k^m||Y_{\xi}||^2,$$
 where $|| . ||$ in the right hand side denotes the $L^2$- norm
 with respect to the K\"ahler form $\widetilde{\omega}$ on $Y$ and
 Fubini-Study metric on the fibres.

\end{prop}

\begin{proof}

By \eqref{eq7}, we have $$|Y_{\xi}|^2=-4\frac{(z^*\xi z)^2- (z^*
\xi^2 z)(z^* z)}{(z^* z)^2}$$ This implies that
\begin{align*}||Y_{\xi}||_{L^2(\widetilde{\omega})}^2&= tr \Big(   \xi^* \xi \int
\frac{zz^*}{z^* z} \widetilde{\omega}^n \Big)-\int \frac{(z^* \xi
z)^2}{(z^* z)^2}\widetilde{\omega}^n\\&=tr \Big(   \xi^* \xi \int
\frac{zz^*}{z^* z} \widetilde{\omega}^n \Big)-\int
\dot{\varphi}^2\widetilde{\omega}^n.\end{align*} We can write
$$\int_{Y} \frac{zz^*}{z^*
z}\widetilde{\omega}^n=D^{(k)}I+M^{(k)},$$ where
$D^{(k)}\rightarrow n_{2}/n_{1}$ as $k\rightarrow \infty$ and
$M^{(k)}$ is a trace free hermitian matrix with
$||M^{(k)}||_{\textrm{op}} \leq \epsilon.$ Therefore,
$$||Y_{\xi}||^2=|\xi|^2D^{(k)}+ tr(\xi^* \xi M^{(k)})-\int
\dot{\varphi}^2\widetilde{\omega}^n.$$ Hence
$$|tr(\xi^* \xi M^{(k)})|=|tr( \xi M^{(k)})\xi| \leq ||\xi||^2
||M^{(k)}||_{\textrm{op}} \leq \epsilon ||\xi||^2.$$ Since
$D^{(k)}\rightarrow n_{2}/n_{1}$ as $k \rightarrow \infty$, there
exists a positive constant $c$ such that $$||Y_{\xi}||^2 \geq
c||\xi||^2-\int \dot{\varphi}^2\widetilde{\omega}^n.$$ On the
other hand \begin{align*}\Big|\int
\dot{\varphi}\widetilde{\omega}^n \Big|= |tr(\xi M^{(k)})| &\leq
\sqrt{N} ||\xi|| ||M^{(k)}||_{\textrm{op}}\\& \leq c
k^{\frac{m}{2}}||\xi|| ||M^{(k)}||_{\textrm{op}}.\end{align*} Now
applying Lemma \ref{lem2}, we get

\begin{align*}C\int_{Y} \dot{\varphi}^2 \widetilde{\omega}^n &\leq k^m\int_{Y}
\overline{\partial}\dot{\varphi}\wedge \partial \dot{\varphi}
\wedge \widetilde{\omega}^{n-1} +k^{-m}\Big( \int_{Y}
\dot{\varphi} \widetilde{\omega}^n \Big)^2\\&\leq k^m\int_{Y}
\overline{\partial}\dot{\varphi}\wedge \partial \dot{\varphi}
\wedge \widetilde{\omega}^{n-1}+c_{2} ||\xi||^2
||M^{(k)}||_{\textrm{op}}^2.\end{align*} This implies
$$(c_{1}-C_{2} ||M^{(k)}||_{\textrm{op}}^2 )||\xi||^2 \leq ||Y_{\xi}||^2 +
k^m  \int_{Y} \overline{\partial}\dot{\varphi}\wedge
\partial \dot{\varphi} \wedge \widetilde{\omega}^{n-1}.$$
Since $||M^{(k)}||_{\textrm{op}} \leq\epsilon$ and $\epsilon $ is
small enough, there exists a positive constant $c$ such that $$c
||\xi||^2 \leq ||Y_{\xi}||^2+ k^m  \int_{Y}
\overline{\partial}\dot{\varphi}\wedge
\partial \dot{\varphi} \wedge \widetilde{\omega}^{n-1}$$
$$= ||Y_{\xi}||^2+ k^m  \int_{Y} |\overline{\partial}\dot{\varphi}|_{\widetilde{\omega}}^2 \widetilde{\omega}^n$$
We know that $
\overline{\partial}\dot{\varphi}|_{Y}=\iota_{\pi_{T}Y_{\xi}}\widetilde{\omega}$
which implies $$c||\xi||^2 \leq ||Y_{\xi}||^2+ k^m
||\pi_{T}Y_{\xi}||^2.$$

\end{proof}

\begin{lem}\label{lem3}

For $k \gg 0 $, we have

$$|S(\omega_{0}+k\omega_{\infty}) |\leq C \log k,$$
where $S$ is the scalar curvature.

\end{lem}

\begin{proof}

We have $$(\omega_{0}+k \omega_{\infty})^n=\sum_{j=0}^{m} {n
\choose k} k^j \omega_{0}^{n-j}\wedge
\omega_{\infty}^j=\big(1+\sum_{j=1}^m k^j f_{j}  \big)
\omega_{0}^n,$$ for some smooth nonnegative functions $f_{j}$ on
$Y$. The function $f_{m}$ is positive, since
$\omega_{0}^{n-m}\wedge \omega_{\infty}^m$ is a volume form.
Therefore there exists a positive constant $l$ such that $f_{m}
\geq l >0$. We define $$F=\sum_{j=1}^m k^{j-m} f_{j}.$$

We have $$\nabla^2
\log(1+k^mF)=\nabla \Big(\frac{k^m \nabla F}{1+k^mF}\Big)=
\frac{k^m\nabla^2F}{1+k^mF}-\frac{k^{2m}(\nabla
F)^2}{(1+k^mF)^2}.$$ Hence there exists a positive constant $C$
such that
$$\big| \log(1+k^mF) \big|_{C^2} \leq mC\log k+ C,$$
since $||F||_{C^2}$ is bounded independent of $k$ and $F \geq f_{m
}\geq l>0$. This implies that
\begin{align*}\big|\partial\overline{\partial} \log\det(\omega_{0}+k
\omega_{\infty})  \big|_{C^0}&\leq \big| \log\det(\omega_{0}+k
\omega_{\infty})  \big|_{C^2} \\&=\big| \log(\omega_{0}+k
\omega_{\infty})^n  \big|_{C^2} \\&\leq \big| \log\omega_{0}^n
\big|_{C^2}+\big| \log (1+k^mF)  \big|_{C^2}\\&\leq  C_{1}+C_{2}m
\log k.\end{align*}

Fix a point $p \in Y$ and a holomorphic local coordinate
$z_{1},...,z_{n}$ around $p$ such that $$\omega_{0}(p)=i \sum
dz_{i} \wedge d\overline{z_{i}},$$
$$\omega_{\infty}(p)=i \sum \lambda_{i}dz_{i} \wedge d\overline{z_{i}},$$
where $\lambda_{i}$'s are some nonnegative real numbers.
Therefore, we have \begin{align*}\big|S(\omega_{0}+k
\omega_{\infty})(p) \big|&=\big| \sum
(1+k\lambda_{i})^{-1}\partial_{i}\partial_{\overline{i}}
\log\det(\omega_{0}+k \omega_{\infty})  \big|\\&\leq \sum
(1+k\lambda_{i})^{-1} (C_{1}+C_{2}m \log k)\leq C_{3} \log k
,\end{align*} for $k \gg 0.$

\end{proof}

\begin{prop}\label{prop4}
 For any holomorphic vector field $V$ on $\mathbb{P}^{N-1}$, we
 have
 $$|\pi_{\mathcal{N}} V|^2 \geq C_{R} k^{-1} |\overline{\partial}(\pi_{\mathcal{N}}
 V)|^2.$$

\end{prop}

\begin{proof}
The following is from (\cite[pp. 705-708]{PS2}). For the sake of
completeness, we give the details of the proof. Fix $x \in Y$. Let
$e_{1},...,e_{n},f_{1},...,f_{m}$ be a local holomorphic frame for
$\iota^* T\mathbb{P}^{N-1}$ around $x$ such that
\begin{enumerate}
\item $e_{1}(x),...,e_{n}(x),f_{1}(x),...,f_{m}(x)$ form an
orthonormal basis. \item $e_{1},...,e_{n}$ is a local holomorphic
basis for $TY$.
\end{enumerate}
Then there exist holomorphic functions $a_{j}$ and $b_{j}$'s such
that $$V=\sum a_{j}e_{j}+\sum b_{j}f_{j}.$$ Notice that
$\pi_{\mathcal{N}}f_{j}-f_{j}$ is tangent to $Y$ since
$\pi_{\mathcal{N}}(\pi_{\mathcal{N}}f_{j}-f_{j})=0.$ Therefore, we
can write $$\pi_{\mathcal{N}}f_{j}-f_{j}=\sum \phi_{ij}e_j,$$
where $\phi_{ij}$'s are smooth functions. Since
$e_{1}(x),...,e_{n}(x),f_{1}(x),...,f_{m}(x)$ form an orthonormal
basis, we have $\phi_{ij}(x)=0$. Then $$\pi_{\mathcal{N}} V=
\sum_{j=1}^{m}b_{j} \big( f_{j}-\sum_{i} \phi_{ij} e_{i}  \big).$$
It implies that $$\overline{\partial}(\pi_{\mathcal{N}} V)=
\sum_{j=1}^{m}b_{j} \big(-\sum_{i} (\overline{\partial}\phi_{ij})
e_{i} \big).$$ So in order to establish \ref{prop4}, we need to
prove that
$$\sum_{i=1}^{n}|\sum_{j=1}^{m}b_{j}\overline{\partial}\phi_{ij}|^2 \leq C_{R}^{-1}k \sum_{j=1}^{m} |b_{j}|^2.$$
Using the Cauchy-Schwartz inequality, it suffices to prove
$$\sum_{i=1}^{n}\sum_{j=1}^{m}|\overline{\partial}\phi_{ij}|^2
\leq C_{2}k,$$ where $C_{2}=C_{2}(R)$ is independent of $k$
(depends on $R$.) Now the matrix
$A^*=(\overline{\partial}\phi_{ij})$ is the dual of the second
fundamental form $A$ of $TY$ in $\iota^* T\mathbb{P}^{N-1}$. Let
$F_{\iota^* T\mathbb{P}^{N-1}}$ be the curvature tensor of the
bundle $\iota^* T\mathbb{P}^{N-1}$ with respect to the
Fubini-Study metric. $\displaystyle F_{\iota^* T\mathbb{P}^{N-1}}$
is a $2$-form on $Y$ withe values in $End(\iota^*
T\mathbb{P}^{N-1})$. Thus $F_{\iota^* T\mathbb{P}^{N-1}}\big
|_{TY}$ is a two form on $Y$ with values in $Hom(TY,\iota^*
T\mathbb{P}^{N-1})$. So, $\pi_{T} \circ (F_{\iota^*
T\mathbb{P}^{N-1}} \big|_{TY})$ is a two form on $Y$ with values
in $End(TY).$ Also let $F_{TY}$ be the curvature tensor of the
bundle $TY$ with respect to the pulled back Fubini-Study metric
$\widetilde{\omega}=\iota^* \omega_{\textrm{FS}}$. Now by
computations in \cite[5.28]{PS2}, we have
$$\sum_{i=1}^{n}\sum_{j=1}^{m}|\overline{\partial}\phi_{ij}|^2= \Lambda_{\widetilde{\omega}} Tr \big[   \pi_{T} \circ (F_{\iota^*
T\mathbb{P}^{N-1}} \big|_{TY})-F_{TY}            \big],$$ where
$\Lambda_{\widetilde{\omega}}$ is the contraction with the
K\"ahler form $\widetilde{\omega}$. The formula \cite[5.33]{PS2}
gives
$$\Lambda_{\widetilde{\omega}} Tr \big[   \pi_{T} \circ (F_{\iota^*
T\mathbb{P}^{N-1}} \big|_{TY})            \big]=n+1.$$ On the
other hand $\Lambda_{\widetilde{\omega}} Tr (   F_{TY} )$ is the
scalar curvature of the metric $\widetilde{\omega}$ on $Y$. Since
$\widetilde{\omega}$ has $R$-bounded geometry, we have
$$|S(\widetilde{\omega})-S(\widetilde{\omega}_{0})| \leq R.$$
Lemma \ref{lem3} implies that $|S(\widetilde{\omega}_{0})| \leq
C\log k \leq Ck$.

\end{proof}

The only thing we need in addition is the following

\begin{prop}\label{prop5}

Assume that there are no nonzero holomorphic vector fields on $Y$.
Then there exists a constant $c^{'}_{R}$ such that for any $\xi
\in su(N)$, we have
$$ c^{'}_{R}||\pi_{T}Y_{\xi}||^2\leq k^{m+2}||\pi_{\mathcal{N}}Y_{\xi}||^2.$$

\end{prop}

\begin{proof}

We define $\alpha=\omega_{0}+\omega_{\infty}$. Since there are no
holomorphic vector fields on $Y$, for any smooth smooth vector
field $W$ on $Y$, we have
$$c||W||^2_{L^2(\alpha)} \leq ||\overline{\partial}W||^2_{L^2(\alpha )}.$$

The trivial inequalities $k\alpha \geq \widetilde{\omega}_{0}$ and
$k^{-m}\widetilde{\omega}_{0}^n \leq \alpha^n \leq
\widetilde{\omega}_{0}^n $ imply that
\begin{align*}c||W||^2_{L^2(\widetilde{\omega}_{0})}=c \int
|W|_{\widetilde{\omega}_{0}}^2 \widetilde{\omega}_{0}^n& \leq
ck^{m+1}\int |W|_{\alpha}^2 \alpha^n \\&\leq k^{m+1}\int
|\overline{\partial}W|_{\alpha}^2 \alpha^n\\&\leq k^{m+1}\int
|\overline{\partial}W|_{\widetilde{\omega}_{0}}^2
\widetilde{\omega}_{0}^n\\&=k^{m+1}
||\overline{\partial}W||^2_{L^2(\widetilde{\omega}_{0})}.\end{align*}
Hence, there exists a positive constant $c$ depends on $R$ and
independent of $k$, such that for any $\widetilde{\omega}_{0}$
having $R$-bounded geometry, we have

$$c||W||^2_{L^2(\widetilde{\omega})} \leq k^{m+1}
||\overline{\partial}W||^2_{L^2(\widetilde{\omega})}.$$
 Now, putting
$W=\pi_{T}Y_{\xi}$, we get

$$c||\pi_{\mathcal{T}}Y_{\xi}||^2_{L^2(\widetilde{\omega})}\leq
k^2||\overline{\partial}(\pi_{\mathcal{T}}Y_{\xi})||^2_{L^2(\widetilde{\omega})}.$$
On the other hand $$||\pi_{\mathcal{N}} V||^2 \geq C_{R} k^{-1}
||\overline{\partial}(\pi_{\mathcal{N}} V)||^2,$$ which implies
the desired inequality.

\end{proof}

\section{Perturbing To A Balanced Metric}
We continue with the notation of the previous section. The goal of
this section is to prove Theorem \ref{newthm2} which gives a
condition for when an almost balanced metric can be perturbed to a
balanced one. In order to do this, first we need to establish
Theorem \ref{newthm1}.  We need the following estimate.

\begin{prop}\label{newprop1}
There exist positive real numbers  $K_{j}$ depends only on
$h_{0}$, $g_{\infty}$ and $j$ such that for any $s \in
H^{0}(Y,\mathcal{O}(1)\otimes L^{k})$, we have $$|\nabla^j
s|_{C^0(\widetilde{\omega}_{0})}^2 \leq K_{j} k^{n+j} \int_{Y}
|s|^2 \frac{\omega_{0}^n}{n!}.$$

\end{prop}

In order to prove Proposition \ref{newprop1}, we start with some
complex analysis.

Let $\varphi$ be a strictly plurisubharmonic function and $\psi$
be a  plurisubharmonic function on $B=B(2) \subset \mathbb{C}^n$
such that $\varphi(0)=\psi(0)=0$. We can find a coordinate on
$B(2)$ such that $$\varphi(z)= |z|^2+O(|z|^2) \,\,\,  \textrm{
and} \,\,\, \psi(z)= \sum \lambda_{i}|z_{i}|^2+O(|z|^2),$$ where
$\lambda_{i} \geq 0$. For any function $u: B\rightarrow
\mathbb{C}$, we define $u^{(k)}(z)=u(\frac{z}{\sqrt{k}}) $.

\begin{thm}
(Cauchy Estimate cf.\cite[Theorem 2.2.3]{Ho})There exist positive
real numbers $C_{j}$ such that for any holomorphic function
$u:B\rightarrow \mathbb{C}$, we have

$$|\nabla^ju|^2(0)\leq C_{j} \int_{|z|\leq1}|u(z)|^2 dz\wedge d\overline{z}$$

\end{thm}

\begin{thm}\label{thm5}
There exist positive real numbers $c_{j}$ depends only on $j,
\varphi, \psi $ and $d\mu$ such that for any holomorphic function
$u:B\rightarrow \mathbb{C}$, we have
$$|\nabla^ju|^2(0)\leq c_{j}k^{n+j} \int_{B(1)}|u|^2 e^{-\varphi-k\psi}d\mu,$$
where $d\mu$ is a fixed volume form on $B$.

\end{thm}

\begin{proof}

Applying Cauchy estimate to $u^{(k)}$, we get
\begin{align*}k^{-j}|\nabla^ju|^2(0)&\leq C_{j} \int_{|z|\leq1}|u^{(k)}(z)|^2
dz\wedge d\overline{z}\\& \leq C \int_{|z|\leq1}|u^{(k)}(z)|^2
e^{-\sum (\lambda_{i}+1)|z_{i}|^2} dz\wedge
d\overline{z},\end{align*} since $e^{-\sum
(\lambda_{i}+1)|z_{i}|^2}$ is bounded from below by a positive
constant on the unit ball. Using the change of variable
$w=\frac{z}{\sqrt{k}}$ we get

\begin{align*}k^{-j}|\nabla^ju|^2(0)&\leq C k^n \int_{|w|\leq k^{-1/2}} |u(w)|^2
e^{-k\sum (\lambda_{i}+1)|w_{i}|^2} dw\wedge d\overline{w}\\&\leq
C k^n \int_{|w|\leq k^{-1/2}} |u(w)|^2 e^{-\sum
(k\lambda_{i}+1)|w_{i}|^2} dw\wedge d\overline{w}.\end{align*} On
the other hand, we have $$\varphi(z)+k\psi(z)= k\sum
(\lambda_{i}+1)|z_{i}|^2+ \mu(z)+k\sigma(z),$$ where
$\displaystyle\mathop {\lim }\limits_{z \to
0}\frac{\mu(z)}{|z|^2}=\displaystyle\mathop {\lim }\limits_{z \to
0}\frac{\sigma(z)}{|z|^2}=0.$

Let $|w|\leq k^{-1/2} $, we have $$|k\sigma(w)+\mu(w) | \leq c
(k|w|^2+|w|^2)\leq 2c$$ for some constant $c$ depending only on
$\psi$ and $\varphi$. Hence

\begin{align*}k^{-j}&|\nabla^ju|^2(0)\leq C k^n \int_{|w|\leq k^{-1/2}}
|u(w)|^2 e^{-\sum (k\lambda_{i}+1)|w_{i}|^2} dw\wedge
d\overline{w} \\
&= Ce^{2c}k^n \int_{|w|\leq k^{-1/2}} |u(w)|^2 e^{-\sum
(k\lambda_{i}+1)|w_{i}|^2-2c}dw\wedge d\overline{w} \\
&\leq C^{'}k^n\int_{|w|\leq k^{-1/2}} |u(w)|^2 e^{-\sum
(k\lambda_{i}+1)|w_{i}|^2-(\mu(w)+k\sigma(w))}dw\wedge
d\overline{w} \\
&=C^{'}k^n \int_{|w|\leq k^{-1/2} } |u(w)|^2
e^{-(\varphi(w)+k\psi(w))}dw\wedge d\overline{w} \\
&\leq C^{'}k^{n}
\int_{B(1)}|u|^2 e^{-\varphi-k\psi}dz\wedge
d\overline{z}.\end{align*}

Hence,

$$|\nabla^ju|^2(0)\leq c_{j} k^{n+j}\int_{B(1)}|u|^2 e^{-\varphi-k\psi}d\mu.$$

\end{proof}

\begin{proof}[Proof of Proposition \ref{newprop1}]
Fix a point $p$ in $Y$ and a geodesic ball $B \subset Y$ centered
at $p$. Let $e_{L}$ be a holomorphic frame for $L$ on $B$ and $e$
be a holomorphic frame for $\mathcal{O}(1)$ such that
$||e_{L}||(p)=||e||(p)=1.$ Any $s \in
H^{0}(Y,\mathcal{O}(1)\otimes L^{k})$ can be written as $s=u e
\otimes e_{L}^{\otimes k}$ for some holomorphic function $u:
B\rightarrow \mathbb{C}$. We have $$\nabla^j s= \sum {j\choose i}
\nabla^iu \otimes \nabla^{j-i}(e \otimes e_{L}^{\otimes k}).$$
Therefore,
$$|\nabla^j s|^2(p)\leq C(\sum |\nabla^iu|^2(p) ||\nabla^{j-i}(e
\otimes e_{L}^{\otimes k})||^2(p) ).$$ On the other hand we have
$$||\nabla^{\alpha}(e \otimes e_{L}^{\otimes k})||^2(p) \leq
\sum_{i=0}^{\alpha}(||\nabla^i e||^2(p)+ k^{\alpha-i}
||\nabla^{\alpha-i} e_{L}||^2(p)) \leq C_{\alpha}k^{\alpha}. $$
Hence $$|\nabla^j s|^2(p)\leq C'(\sum |\nabla^iu|^2(p) k^{j-i}
).$$ Applying Theorem \ref{thm5} concluds the proof.

\end{proof}

For the rest of this section, we fix a positive integer $q$. We
continue with the notation ($Y, \omega_{\infty}, \omega_{0},
\widetilde{\omega_{0}}$ ) of section $3$. In the rest of this
section, we fix the reference metric $\omega_{0}$ on $Y$ and
recall the Definition \ref{def1}.

\begin{Def}\label{newdef1}
The sequence of hermitian metrics $h_{k} $ on
$\mathcal{O}(1)\otimes L^{k}$ and ordered bases
$\underline{s}^{(k)}=(s_{1}^{(k)},...,s_{N}^{(k)})$ for
$H^0(Y,\mathcal{O}(1)\otimes L^{k})$ is called \emph{almost
balanced of order $q$} if for any $k$  $$\sum
|s_{i}^{(k)}|_{h_{k}}^2=1$$ and
$$\int_{Y}\langle s_{i}^{(k)}, s_{j}^{(k)} \rangle
{h_{k}}dvol_{h_{k}}=D^{(k)}\delta_{ij}+M^{(k)}_{ij},$$ where
$D^{(k)}$ is a scalar so that $D^{(k)} \rightarrow n_{2}/n_{1}$ as
$k \rightarrow \infty$ (See \eqref{eq4} and \eqref{eq14}.), and
$M^{(k)}$ is a trace-free hermitian matrix such that
$||M^{(k)}||_{\textrm{op}}=O(k^{-q-1})$.

\end{Def}
We state the following lemma without proof. The proof is a
straightforward calculation.

\begin{lem}\label{lem4}

Let the sequence of hermitian metrics $h_{k} $ on
$\mathcal{O}(1)\otimes L^{k}$ and ordered bases
$\underline{s}^{(k)}=(s_{1}^{(k)},...,s_{N}^{(k)})$ for
$H^0(Y,\mathcal{O}(1)\otimes L^{k})$ be almost balanced of order
$q$. Suppose
\begin{equation}\label{neq1}||\widetilde{\omega_{k}}-\widetilde{\omega_{0}}||_{C^a(\widetilde{\omega_{0}})}=O(k^{-1}),\end{equation}
where $\widetilde{\omega_{k}}= \textrm{Ric}(h_{k}).$ Then for any
$\epsilon >0$ there exists a positive integer $k_{0}$ such that
$$\widetilde{\omega_{k}} \geq (1-\epsilon) \widetilde{\omega_{0}} \,\,\,  \textrm{for} \,\,\, k \geq k_{0}.$$

\end{lem}

Assume that there exist a sequence of almost balanced metrics
$h_{k}$ of order $q$  and bases
$\underline{s}^k=(s_{1}^{(k)},...,s_{N}^{(k)})$ for
$H^{0}(Y,\mathcal{O}(1)\otimes L^{k})$ which satisfies
\eqref{neq1}. As before $\widetilde{\omega_{k}}=
\textrm{Ric}(h_{k}).$  Then Lemma \ref{lem4} implies that for $k
\gg 0$, $\widetilde{\omega_{k}}$ has $R$-bounded geometry.

Fix $k$ and let $B \in isu(N_{k})$. Without loss of generality, we
can assume that $B$ is the diagonal matrix
$\textrm{diag}(\lambda_{i}),$ where $\lambda_{i} \in \mathbb{R}$
and $\sum \lambda_{i}=0$. There exists a unique hermitian metric
$h_{B}$ on $\mathcal{O}_{\mathbb{P}E^*}(1)\otimes L^{k}$ such that
$$\sum e^{2\lambda_{i}}|s_{i}^{(k)}|_{h_{B}}^2=1.$$ Let
$\widetilde{\omega}_{B}=Ric (h_{B})$. In the next theorem, we will
prove that there exist a constant $c$ and open balls $U_{k}
\subset isu(N_{k})$ around the origin of radius $ck^{-(n+a+2)}$ so
that if $B \in U_{k}$, then $h_{B}$ is $R$-bounded. More
precisely,
\begin{thm}\label{newthm1}
Suppose that \eqref{neq1} holds.
\begin{itemize}

\item There exist $c >0$ and $k_{0} >0$ such that if $k \geq
k_{0}$ and $B \in isu(N_{k})$ satisfies $$||B||_{\textrm{op}} \leq
ck^{-(n+a+2)}R,$$ then the metric $\widetilde{\omega_{B}}$ is
$R$-bounded.

\item There exists $c >0$ such that if $B \in isu(N_{k})$
satisfies $$||B||_{\textrm{op}} \leq k^{-(n+a+3)},$$ then
$$||M^{B}||_{\textrm{op}} \leq ck^{-1},$$ where the matrix $M^{B}=
(M^B_{ij})$is defined by
$$M^B_{ij}=e^{\lambda_{i}+\lambda_{j}}\int_{Y}\langle
s_{i}^{(k)},s_{j}^{(k)}\rangle_{h_{B}}\frac{\widetilde{\omega}_{B}^n}{n!}-\frac{V_{k}}{N_{k}}\delta_{ij}.$$

\end{itemize}

\end{thm}

\begin{proof}

Let $h_{B}=e^{\varphi_{B}}h_{k}$. So, we have $$1=\sum
e^{2\lambda_{i}}|s_{i}^{(k)}|_{h_{B}}^2=e^{\varphi_{B}}\sum
e^{2\lambda_{i}}|s_{i}^{(k)}|_{h_{k}}^2.$$ Hence
$$\varphi_{B}=-\log \sum e^{2\lambda_{i}}|s_{i}^{(k)}|_{h_{k}}^2=-\log
\big(1+ \sum (e^{2\lambda_{i}}-1)|s_{i}^{(k)}|_{h_{k}}^2\big).$$



If $||B||_{\textrm{op}}$ is small enough, there exists $C >0 $ so
that
$$||\varphi_{B}||_{C^{a+2}(\widetilde{\omega}_{0})}\leq C ||B||_{\textrm{op}}\sum_{i=1}^{N_{k}}|\nabla^{a+2}  s_{i}^{(k)}|_{C^0(\widetilde{\omega}_{0})}^2$$
and therefore Proposition \ref{newprop1} implies that
\begin{align*}||\varphi_{B}||_{C^{a+2}(\widetilde{\omega}_{0})}&\leq C
||B||_{\textrm{op}}k^{n+a+2}\sum_{i=1}^{N_{k}}
\int_{Y}|s_{i}^{(k)}|_{h_{k}}^2\frac{\omega_{0}^n}{n!}\\&=C
||B||_{\textrm{op}}k^{n+a+2}
\int_{Y}\sum_{i=1}^{N_{k}}|s_{i}^{(k)}|_{h_{k}}^2\frac{\omega_{0}^n}{n!}\\&
=C||B||_{\textrm{op}}k^{n+a+2}\int_{Y}\frac{\omega_{0}^n}{n!}=c_{1}||B||_{\textrm{op}}k^{n+a+2}\end{align*}
for some positive constant $c_{1}$. Now if $||B||_{\textrm{op}}
\leq c_{1}^{-1}\frac{R-1}{2R}k^{-(n+a+2)}$, then
\begin{equation}\label{neq2}  ||\varphi_{B}||_{C^{a+2}(\widetilde{\omega}_{0})}\leq \frac{R-1}{2R}.  \end{equation}
Therefore, $$||i\overline{\partial}\partial
\varphi_{B}||_{C^{0}(\widetilde{\omega}_{0})}\leq
\frac{R-1}{2R},$$ which implies that
\begin{equation} \label{neq3}  i\overline{\partial}\partial
\varphi_{B} \geq  -\frac{R-1}{2R}\widetilde{\omega_{0}} .
\end{equation}

In order to show that $\widetilde{\omega_{B}}$ is $R$-bounded, we
need to prove the following: \begin{align}\label{neq4}
&||\widetilde{\omega}-\widetilde{\omega}_{0}||_{C^{a}(\widetilde{\omega}_{0})}\leq
R,\\& \label{neq5}\widetilde{\omega_{B}} \geq \frac{1}{R}
\widetilde{\omega}_{0}.\end{align} To prove \eqref{neq4},
\eqref{neq1} and \eqref{neq2} imply that for $k \gg 0$

\begin{align*}||\widetilde{\omega}_{B}-\widetilde{\omega}_{0}||_{C^a(\widetilde{\omega}_{0})}&\leq
||\widetilde{\omega}_{B}-\widetilde{\omega}_{k}||_{C^a(\widetilde{\omega}_{0})}+||\widetilde{\omega}_{k}-\widetilde{\omega}_{0}||_{C^a(\widetilde{\omega}_{0})}\\&\leq
||\varphi_{B}||_{C^{a+2}(\widetilde{\omega}_{0})}+k^{-1} \leq
\frac{R-1}{2R}+k^{-1}\\& \leq R.\end{align*}

To prove \eqref{neq5}, applying Lemma \ref{lem4} with
$\epsilon=\frac{R-1}{2R}$ gives $$\widetilde{\omega_{k}} \geq
\frac{R+1}{2R} \widetilde{\omega_{0}},$$ and therefore
\eqref{neq3} implies $$\widetilde{\omega_{B}}-\frac{1}{R}
\widetilde{\omega}_{0}
=\widetilde{\omega_{k}}+i\overline{\partial}\partial \varphi_{B}-
\frac{1}{R} \widetilde{\omega}_{0} \geq
\widetilde{\omega_{k}}-\frac{R+1}{2R}\widetilde{\omega_{0}} \geq
0,$$ for $k \gg 0.$

In order to prove the second part, by a unitary change of basis,
we may assume without loss of generality that the matrix $M^{B}$
is diagonal. By definition

$$M^B_{ij}=e^{\lambda_{i}+\lambda_{j}} \int_{Y} F  \langle
s_{i}, s_{j}\rangle
\frac{\widetilde{\omega}_{B}^n}{n!}-\frac{V_{k}}{N_{k}}\delta_{ij},
$$ where
$$F= e^{-\varphi_{B}}\frac{\widetilde{\omega_{B}}^n}{\widetilde{\omega_{k}}^n}.$$
We have \begin{align*}M^B_{ii}&=e^{2\lambda_{i}}\int_{Y} F
|s_{i}|^2_{h_{k}}\frac{\widetilde{\omega}_{k}^n}{n!}-\frac{V_{k}}{N_{k}}\delta_{ij}\\&=e^{2\lambda_{i}}\int_{Y}F
|s_{i}|^2_{h_{k}}\frac{\widetilde{\omega}_{k}^n}{n!}-\int_{Y}
|s_{i}|^2_{h_{k}}\frac{\widetilde{\omega}_{k}^n}{n!}+(M^{(k)})_{ii}\\&=\int_{Y}(e^{2\lambda_{i}}F-1)
|s_{i}|^2_{h_{k}}\frac{\widetilde{\omega}_{k}^n}{n!}+(M^{(k)})_{ii}.\end{align*}
Therefore,
$$|M^B_{ii}| \leq ||e^{2\lambda_{i}}F-1||_{\infty}(\int_{Y}|s_{i}|^2_{h_{k}}\frac{\widetilde{\omega}_{k}^n}{n!})+|(M^{(k)})_{ii}| \leq C(||e^{2\lambda_{i}}F-1||_{\infty}+k^{-q-1}).$$
Define $f=\displaystyle
\frac{\widetilde{\omega_{B}}^n}{\widetilde{\omega_{k}}^n}.$ If
$||B||_{\textrm{op}} \leq k^{-(n+a+3)}$, then
$$|f-1|=\big|\frac{\widetilde{\omega_{B}}^n-\widetilde{\omega_{k}}^n}{\widetilde{\omega_{k}}^n}\big| =O(k^{-1})$$
and
$$|(e^{2\lambda_{i}-\varphi_{B}}-1)|=O(k^{-1}).$$
Therefore,
\begin{align*}||e^{2\lambda_{i}}F-1||&=||e^{2\lambda_{i}-\varphi_{B}}\frac{\widetilde{\omega_{B}}^n}{\widetilde{\omega_{k}}^n}-1||=||e^{2\lambda_{i}-\varphi_{B}}f-1||\\&\leq||(e^{2\lambda_{i}-\varphi_{B}}-1)(f-1)||+||(f-1)||\\&\qquad +||(e^{2\lambda_{i}-\varphi_{B}}-1)||\\&=O(k^{-1}),\end{align*}
which implies that
$$||M^{B}||_{\textrm{op}}=O(k^{-1}).$$

\end{proof}

\begin{thm}\label{newthm2}

Suppose that the sequence of metrics $h_{k}$ on
$\mathcal{O}(1)\otimes L^{k}$ and bases
$\underline{s}^k=(s_{1}^k,...,s_{N}^k)$ for
$H^{0}(Y,\mathcal{O}(1)\otimes L^{k})$  is almost balanced of
order $q$. Suppose that \eqref{neq1} holds for
$$\widetilde{\omega}_{k}=Ric(h_{k}) \,\,\, \textrm{ and} \,\,\,
\omega_{k}=Ric(h_{k})-k\omega_{\infty}.$$ If $q >
\frac{5m}{2}+n+a+5$, then $(Y,\mathcal{O}(1)\otimes L^{k})$admits
balanced metric for $k \gg 0.$

\end{thm}

\begin{proof}

Let $R >1$ and $k$ be a fixed large positive integer. Let $\sigma
\in isu(N)$, where $N=N_{k}=\dim H^{0}(Y,\mathcal{O}(1)\otimes
L^{k}).$ If $||\sigma||_{\textrm{op}} \leq
\frac{c}{2}k^{-(n+a+3)}R$, then Theorem \ref{newthm1} implies that
$e^{\sigma}\underline{s}$ has $R$-bounded geometry and
$||M^{\sigma}||_{\textrm{op}}\leq \epsilon$ for $k \gg 0$, where
$\epsilon $ is the constant in the statement of Theorem
\ref{thm4}. Thus, Theorem \ref{thm4} implies that
$\Lambda(e^{\sigma}\underline{s}^{(k)})\leq C k^{2m+2}=\lambda.$
With the notation of Proposition \ref{prop1}, we have
$\mu(z_{0})=M^{(k)}$. Therefore $$|\mu(z_{0})|=|M^{(k)}| \leq
\sqrt{N_{k}}||M^{(k)}||_{\textrm{op}}\leq
C^{'}k^{\frac{m}{2}-q}.$$

Letting $\delta= \frac{c}{2} k^{-(n+a+3)}R$, we have $\lambda
|\mu(z_{0})| <\delta $ if $q > \frac{5m}{2}+n+a+5 $ and $k \gg 0$.
Therefore if $q > \frac{5m}{2}+n+a+5 $ and $k \gg 0$, we can apply
Proposition \ref{prop1} to get balanced metrics for $k \gg 0$.

\end{proof}

\begin{rem}
By Proposition \ref{prop1}, there exists $\sigma_{0}$ such that
$e^{\sigma_{0}}\underline{s}$ is balanced and $|\sigma_{0}|\leq (C
k^{2m+2})(C^{'}k^{\frac{m}{2}-q})=C^{''}k^{\frac{5m}{2}+2-q}.$
Since $$||\varphi_{\sigma_{0}}||_{C^{a+2}}\leq c_{1}k^{n+a+2}
||\sigma_{0}||_{\textrm{op}}\leq c_{1}k^{n+a+2} |\sigma_{0}|,$$
then
$$||\varphi_{\sigma_{0}}||_{C^{a+2}}\leq
ck^{\frac{5m}{2}+n+a+4-q}.$$ Therefore,
\begin{equation}\label{eq15}||\omega_{k}^{\textrm{bal}}-\widetilde{\omega}_{k}||_{C^{a}(\omega_{0})}=O(k^{-1}).\end{equation}
\end{rem}

\section{Asymptotic Expansion}

The goal of this section is to prove Theorem \ref{thmH1}. Theorem
\ref{thmH1} gives an asymptotic expansion for the Bergman kernel
of $(\mathbb{P}E^*,\mathcal{O}_{\mathbb{P}E^*}(1)\otimes
\pi^*L^{k})$. We obtain such an expansion by using the Bergman
kernel asymptotic expansion proved in (\cite{C}, \cite{Z}). Also
we compute the first nontrivial coefficient of the expansion. In
the next section, we use this to construct sequence of almost
balanced metrics. We start with some linear algebra.

Let $V$ be a hermitian vector space of dimension $r$. The
projective space $\mathbb{P}V^*$ can be identified with the space
of hyperplanes in $V$ via $$f\in V \rightarrow
{ker(f)}=V_{f}\subseteq V.$$ If $f\neq 0 $ then $V_{f}$ will be a
hyperplane. There is a natural isomorphism between $V$ and
$H^{0}(\mathbb{P}V^*,\mathcal{O}_{\mathbb{P}V^*}(1))$ which sends
$v\in V$ to $\hat{v}\in
H^{0}(\mathbb{P}V^*,\mathcal{O}_{\mathbb{P}V^*}(1)) $ such that
for any $f \in V^*, \hat{v}(f)=f(v) $. Now we can see that the
inner product on $V$ induces an inner product on $V^*$ and then a
metric on $\mathcal{O}_{\mathbb{P}V^*}(1).$ For $v,w \in V$ and $f
\in V^*$ we define $$<\hat{v},\hat{w}>_{[f]}=
\frac{f(v)\overline{f(w)}}{|f|^{2}}.$$

 \begin{Def}\label{defH1}
For any inner product  $h$ on $V$, We denote the induced metric on
$\mathcal{O}_{\mathbb{P}V^*}(1)$ by $\widehat{h}$.

 \end{Def}

The following is a straight forward computation.

\begin{prop}\label{propH1}
 For any $v,w \in V$ we have
$$<v,w>_{h}=C_{r}^{-1}\int_{\mathbb{P}V^*} <\hat{v},\hat{w}>_{\widehat{h}}
\frac{\omega_{\textrm{FS}}^{r-1}}{(r-1)!}$$ where $C_{r}$ is a
constant defined by
\begin{equation}\label{eq16}C_{r}=\int_{\mathbb{C}^{r-1}}
\frac{d\xi \wedge d\overline{\xi}}{(1+\sum_{j=1}^{r-1}
|\xi_{j}|^{2})^{r+1}}.\end{equation} Here $d\xi \wedge
d\overline{\xi}= (\sqrt{-1}d\xi_{1} \wedge
d\overline{\xi}_{1})\wedge \dots \wedge (\sqrt{-1}d\xi_{r-1}
\wedge d\overline{\xi}_{r-1}).$

\end{prop}

\begin{Def}\label{defH2}
For any $v \in V$, we define an endomorphism of $V$ by
$$\lambda(v,h)=\frac{1}{||v||_{h}^2} v \otimes v^{*_{h}},$$
where $v^{*_{h}}(.)= h(.,v).$

\end{Def}

Let $(X,\omega)$ be a K\"ahler manifold of dimension $m$ and $E$
be a holomorphic vector bundle on $X$ of rank $r$. Let $L$ be an
ample line bundle on $X$ endowed with a Hermitian metric $\sigma$
such that $Ric (\sigma)=\omega$. For any hermitian metric $h$ on
$E$, we define the volume form
$$d\mu_{g}=\frac{\omega_{g}^{r-1}}{(r-1)!}\wedge \frac{\pi^*
\omega^m}{m!},$$ where $g=\widehat{h}$
,$\omega_{g}=Ric(g)=Ric(\widehat{h})$ and
$\pi:\mathbb{P}E^*\rightarrow X$ is the projection map. The goal
is to find an asymptotic expansion for the Bergman kernel of
$\mathcal{O}_{\mathbb{P}E^*}(1) \otimes L^k \rightarrow
\mathbb{P}E^*$ with respect to the $L^2$-metric defined on
$H^{0}(\mathbb{P}E^*,\mathcal{O}_{\mathbb{P}E^*}(1)\otimes
\pi^*L^{k})$. We define the $L^2$- metric using the fibre metric
$g \otimes \sigma^{\otimes k}$ and the volume form $d\mu_{g,k}$
defined as follows
\begin{equation}\label{eqH1}d\mu_{g,k}=k^{-m}\frac{(\omega_{g} + k
\omega)^{m+r-1}}{(m+r-1)!}= \sum_{j=0}^{m} k^{j-m}
\frac{\omega_{g}^{m+r-1-j}}{(m+r-j)!} \wedge
\frac{\omega^j}{j!}.\end{equation}

In order to do that, we reduce the problem to the problem of
Bergman kernel asymptotics on $E\otimes L^k \rightarrow X$. The
first step is to use the volume form $d\mu_{g}$ which is a product
volume form instead of the more complicated one $d\mu_{g,k}$. So,
we replace the volume form $d\mu_{g,k}$ with $d\mu_{g}$ and the
fibre metric $g \otimes \sigma^k$ with $g(k) \otimes \sigma^k$,
where the metrics $g(k)$ are defined on
$\mathcal{O}_{\mathbb{P}E^*}(1)$ by
\begin{equation}\label{eqH2}g(k)=k^{-m}(\sum_{j=0}^{m} k^j
f_{j})g=(f_{m}+k^{-1}f_{m-1}+...+k^{-m}f_{0})g,\end{equation} and
\begin{equation} \label{eqH3}\frac{\omega_{g}^{m+r-1-j}}{(m+r-j)!} \wedge
\frac{\omega^j}{j!}=f_{j}d\mu_{g}. \end{equation} Clearly the
$L^2$-inner products $L^2(g \otimes \sigma^k,d\mu_{g,k})$ and
$L^2(g(k) \otimes \sigma^k,d\mu_{g})$ on
$H^{0}(\mathbb{P}E^*,\mathcal{O}_{\mathbb{P}E^*}(1)\otimes
\pi^*L^{k})$ are the same. The second step is going from
$\mathcal{O}_{\mathbb{P}E^*}(1) \rightarrow \mathbb{P}E^*$ to
$E\rightarrow X$. In order to do this we somehow push forward the
metric $g(k)$ to get a metric $\widetilde{g}(k)$ on $E$ (See
Definition \ref{defH4}). Then we can apply the result on the
asymptotics of the Bergman kernel on $E$. The last step is to use
this to get the result.

\begin{Def}\label{defH3}
Let $\widehat{s_{1}^k},....,\widehat{s_{N}^k} $ be an orthonormal
basis for
$H^{0}(\mathbb{P}E^*,\mathcal{O}_{\mathbb{P}E^*}(1)\otimes
\pi^*L^{k})$ w.r.t. $L^2(g \otimes \sigma^k, d\mu_{k,g})$. We
define
\begin{equation}\label{eqH4}\rho_{k}(g,\omega)=\sum_{i=1}^{N}
|\widehat{s_{i}^k}|_{g\otimes \sigma^k}^2.\end{equation}
\end{Def}

\begin{Def}
For any $(j,j)$-form $\alpha$ on $X$, we define the contraction
$\Lambda_{\omega}^j\alpha$ of $\alpha$ with respect to the
K\"ahler form $\omega$ by $$\frac{m!}{(m-j)!} \alpha \wedge
\omega^{m-j}=(\Lambda_{\omega}^j\alpha )\,\,\omega^m.$$
\end{Def}

In this section we fix the K\"ahler form $\omega$ on $X$ and
therefore simply denote $\Lambda_{\omega}^j\alpha$ by $\Lambda
^j\alpha$.
\begin{lem}\label{lemH5}
Let $\nu_{0}$ be a fixed K\"ahler form on $X$. For any positive
integer $p$ there exists a constant $C$ such that for any
$(j,j)$-form $\gamma$, we have

$$||\nabla^p( \Lambda^j
\gamma)|| \leq \frac{C}{\inf_{x \in
X}|\omega(x)^m|_{\nu_{0}(x)}}(||\gamma||_{C^p(\nu_{0})}+
||\Lambda^j
\gamma||_{C^{p-1}(\nu_{0})})(\sum_{i=1}^{m}||\omega||_{C^p(\nu_{0})}^i).$$

\end{lem}

\begin{proof}

Let $\gamma$ be a $(j,j)$-form. By definition, we have
$$(\Lambda^j \gamma) \,\,\omega^m= \frac{m!}{(m-j)!}\gamma \wedge \omega^{m-j}.$$
Therefore for any positive integer $p$, we have
$$\nabla^p((\Lambda^j \gamma)
\,\,\omega^m)=\frac{m!}{(m-j)!}\nabla^p( \gamma \wedge
\omega^{m-j}).$$ Applying Leibnitz rule, we get
$$\sum_{i=0}^{p} {p \choose i} \nabla^i (\Lambda^j \gamma) \nabla^{p-i}\omega^m  = \frac{m!}{(m-j)!}  \sum_{i=0}^{p} {p \choose i}\nabla^i \gamma \wedge \nabla^{p-i} \omega^{m-i}.$$
Thus there exists a positive constant $C'$ so that
\begin{align*}||\nabla^p( \Lambda^j \gamma)\omega^m||_{C^0(\nu_{0})}\leq
C'( ||\omega^m||_{C^p(\nu_{0})}|| \Lambda^j
\gamma||_{C^{p-1}(\nu_{0})}+||\gamma||_{C^p(\nu_{0})}||\omega^{m-j}||_{C^p(\nu_{0})}).\end{align*}
On the other hand there exists constant $c_{p,j}$ such that for
any any $0 \leq j \leq m-1$,
$$||\omega^{m-j}||_{C^p(\nu_{0})}\leq
c_{p,j}||\omega||_{C^p(\nu_{0})}^{m-j}\leq
c_{p,j}(\sum_{i=1}^{m}||\omega||_{C^p(\nu_{0})}^i).$$ Hence there
exists a constant $C$ such that $$||\nabla^p( \Lambda^j \gamma)||
\leq \frac{C}{\inf_{x \in
X}|\omega(x)^m|_{\nu_{0}(x)}}(||\gamma||_{C^p(\nu_{0})}+
||\Lambda^j
\gamma||_{C^{p-1}(\nu_{0})})(\sum_{i=1}^{m}||\omega||_{C^p(\nu_{0})}^i).$$

\end{proof}

\begin{Def}\label{defH4}
For any hermitian form $g$ on $\mathcal{O}_{\mathbb{P}E^*}(1)$, we
define a hermitian form $\widetilde{g}$ on $E$ as follow
\begin{equation}\label{eqH6}\widetilde{g}(s,t)=C_{r}^{-1} \int_{\mathbb{P}E_{x}^*}
g(\widehat{s},\widehat{t}\, )
\frac{\omega_{g}^{r-1}}{(r-1)!},\end{equation} for $s,t \in
E_{x}.$ (See \eqref{eq16} for definition of $C_{r}$.)
\end{Def}
Notice that if $g=\widehat{h}$ for some hermitian metric $h$ on
$E$, Proposition \ref{propH1} implies that $\widetilde{g}=h.$
Define hermitian metrics $\widetilde{g_{j}}$'s on $E$ by
\begin{equation}\label{eqH7}\widetilde{g_{j}}(s,t)=C_{r}^{-1} \int_{\mathbb{P}E_{x}^*}
 f_{j}g(\widehat{s},\widehat{t} ) \frac{\omega_{g}^{r-1}}{(r-1)!},\end{equation} for
$s,t \in E_{x}.$ Also we define $\Psi_{j} \in End(E)$ by
\begin{equation}\label{eqH8}\widetilde{g_{j}}=\Psi_{j} h.\end{equation}

\begin{prop}\label{propH4}
Let $\nu_{0}$ be a fixed K\"ahler form on $X$ as in Lemma
\ref{lemH5}. For any positive numbers $l$ and $l'$ and any
positive integer $p$, there exists a positive number $C_{l,l',p}$
such that if $$||\omega||_{C^p(\nu_{0})}, ||h||_{C^{p+2}(\nu_{0})}
\leq l $$ and $$\inf_{x \in X}|\omega(x)^m|_{\nu_{0}(x)}\geq l',$$
then
$$||\Psi_{i}||_{C^p(\nu_{0})} \leq C_{l,l',p} \,\,\, \textrm{for any} \,\,\, 1 \leq i\leq m.$$

\end{prop}
\begin{proof}
Fix a point $p \in X$. Let $e_{1},...,e_{r}$ be a local
holomorphic frame for $E$ around $p$ such that $$\langle e_{i},
e_{j}\rangle _{h}(p)=\delta_{ij}, \,\,\,\,\,\,\,\, d \langle
e_{i}, e_{j}\rangle _{h}(p)=0 $$ and $$\frac{i}{2\pi
}F_{h}(p)=\left(
\begin{matrix}
 \omega_{1} & 0& \cdots &0 \\
 0 & \omega_{2}& \cdots & 0\\
 \vdots &  & \ddots  &\vdots \\
 0 & 0  & \cdots & \omega_{r} \end{matrix} \right).$$

Let $\lambda_{1},...,\lambda_{r}$ be the homogeneous coordinates
on the fibre. At the fixed point $p$, we have $$\omega_{g}=
\omega_{\textrm{FS},g}+\frac{\sum \omega_{i} |\lambda_{i}|^2}{\sum
|\lambda_{i}|^2}.$$

Therefore,\begin{align*}\omega_{g}^{r+j-1} \wedge
\omega^{m-j}&={r+j-1 \choose r-1}\omega_{\textrm{FS},g}^{r-1}
\wedge \big(\frac{\sum \omega_{i} |\lambda_{i}|^2}{\sum
|\lambda_{i}|^2}\big)^j \wedge \omega^{m-j}.\end{align*}
Definition of $f_{m-j}$ gives

$$f_{m-j}\omega_{g}^{r-1}\wedge
\omega^m={m \choose j}\omega_{g}^{r-1} \wedge \Big(
\big(\frac{\sum  \omega_{i} |\lambda_{i}|^2}{\sum
|\lambda_{i}|^2}\big)^j \wedge \omega^{m-j} \Big)$$
 Hence

$$f_{m-j}\omega_{\textrm{FS},g}^{r-1}\wedge
\omega^m={m \choose j}\omega_{\textrm{FS},g}^{r-1} \wedge \Big(
\big(\frac{\sum  \omega_{i} |\lambda_{i}|^2}{\sum
|\lambda_{i}|^2}\big)^j \wedge \omega^{m-j} \Big).$$ Therefore,
$$\omega_{\textrm{FS},g}^{r-1}\wedge\Big(f_{m-j}
\omega^m-{m \choose j} \big(\frac{\sum \omega_{i}
|\lambda_{i}|^2}{\sum |\lambda_{i}|^2}\big)^j \wedge \omega^{m-j}
\Big)=0,$$ which implies

\begin{align*}f_{m-j} \omega^m&={m \choose j} \big(\frac{\sum \omega_{i}
|\lambda_{i}|^2}{\sum |\lambda_{i}|^2}\big)^j \wedge
\omega^{m-j}\,\,\,\,\,\,\,\,\,\,\,\,\,\,\,\,\,\,\,\,\,\,\,\,\,\,\,\,\,\,\,\,\,\,\,\,\,\,\,\,\,\,\,\,\,\,\,\,\,\,\,\,\,\,\,\,\,\,\,\,\,\,\,\,\,\,\,\,\,\,\,\,\,\,\,\,\,\,\,\,\,\,\,\,\,\,\,\,\,\,\,
\\&={m \choose j} \frac{\sum_{j_{1}+\dots+j_{r}=j} {j\choose
j_{1},\dots,j_{r}}\omega_{1}^{j_{1}}\wedge \dots \wedge
\omega_{r}^{j_{r}} |\lambda_{1}|^{2j_{1}}\dots
|\lambda_{r}|^{2j_{r}}}{(\sum |\lambda_{i}|^2)^j} \wedge
\omega^{m-j}.\end{align*} Simple calculation gives
$$\int_{\mathbb{C}^{r-1}}
\frac{|\lambda_{\alpha}|^2|\lambda_{1}|^{2j_{1}}\dots
|\lambda_{r-1}|^{2j_{r-1}}d\lambda \wedge
d\overline{\lambda}}{(1+\sum_{j=1}^{r-1}
|\lambda_{j}|^{2})^{r+j+1}}=\frac{C_{r}r!j_{1}!\dots
j_{r}!(j_{\alpha}+1)}{ (r+j)!},$$ when $j_{1}+\dots +j_{r}=j$ and
$1\leq \alpha\leq r$. Hence
\begin{align}\label{neq7}\widetilde{g}_{m-j}(e_{\alpha},e_{\alpha})=C_{r}^{-1}\pi_{*}\big(
f_{m-j}g(\widehat{e_{\alpha}},\widehat{e_{\alpha}})
\frac{\omega_{g}^{r-1}}{(r-1)!}
\big)\\=\frac{r!}{(r+j)!}\Lambda^j\big(
\sum_{j_{1}+\dots+j_{r}=j}(j_{\alpha}+1)\omega_{1}^{j_{1}}\wedge
\dots \wedge \omega_{r}^{j_{r}} \big).\notag\end{align} From
theory of symmetric functions, one can see that there exist
polynomials $P_{i}(x_{1},\dots,x_{j})$ of degree $i$ such that
$$\Psi_{m-j}=\Lambda^j \Big( F_{h}^j+P_{1}(c_{1}(h),\dots,c_{j}(h))F_{h}^{j-1}+\dots+ P_{j}(c_{1}(h),\dots,c_{j}(h)) \Big),$$
where $c_{i}(h)$ is the $i$ th chern form of $h$. Since
$||h||_{C^{p+2}(\nu_{0})} \leq l$, there exists a positive
constant $c'$ such that
$$||F_{h}^j+\dots+
P_{j}(c_{1}(h),\dots,c_{j}(h))||_{C^p(\nu_{0})} \leq c'(1+l)^j.$$
Therefore Lemma \ref{lemH5} implies that
$$||\nabla^p \Psi_{m-j}|| \leq \frac{C}{l'}(c'(1+l)^j+ || \Psi_{m-j}||_{C^{p-1}(\nu_{0})})(1+l)^m,$$
since $$\inf_{x \in X}|\omega(x)^m|_{\nu_{0}(x)}\geq l'$$ and
$$\sum_{i=1}^{m}||\omega||_{C^p(\nu_{0})}^i \leq \sum_{i=1}^{m}l^i \leq
(1+l)^m.$$ On the other hand
\begin{align*}||\Psi_{m-i}||_{C^p(\nu_{0})}&= ||\nabla^p
\Psi_{m-j}||+||\Psi_{m-i}||_{C^{p-1}(\nu_{0})}\\&\leq
\frac{C}{l'}(c'(1+l)^j+ ||
\Psi_{m-j}||_{C^{p-1}(\nu_{0})})(1+l)^m+||\Psi_{m-i}||_{C^{p-1}(\nu_{0})}.\end{align*}
Now we can conclude the proof by induction on $p$.

\end{proof}

\begin{lem}\label{lemH2}
We have the following
\begin{enumerate}
\item $\displaystyle \Psi_{m}=I_{E}.$ \item $\displaystyle
\Psi_{m-1}=\frac{i}{2\pi(r+1)}(Tr(\Lambda F_{h})I_{E}+  \Lambda
F_{h})$.

\end{enumerate}

\end{lem}

\begin{proof}
The first part is an immediate consequence of Proposition
\ref{propH1} and the definition of $\Psi_{m}$. For the second
part, we use the notation used in the proof of Proposition
\ref{propH4}. It is easy to see that for $\alpha \neq \beta$, we
get $\widetilde{g_{m-1}}(e_{\alpha}, e_{\beta})=0$. On the other
hand by plugging $j=1$ in \eqref{neq7}, we get

\begin{align*} \widetilde{g_{m-1}}(e_{\alpha},
e_{\alpha})&=\frac{1}{(r+1)}(Tr(\Lambda F)+ \Lambda
\omega_{\alpha}).\end{align*}

\end{proof}

The following lemmas are straightforward.

\begin{lem}\label{lemH3}

$\displaystyle \widetilde{g \otimes \sigma^k}=\widetilde{g}
\otimes \sigma^k.$

\end{lem}

\begin{lem}\label{lemH4} Let $s_{1},..., s_{N}$ be a basis for $H^0(X,E)$. Then
$$\sum |\widehat{s_{i}}([v^*])|^2_{\widehat{h}}= Tr \big( B \lambda(v^*,h)          \big),$$
where $B=\sum s_{i} \otimes s_{i}^{*_{h}}.$

\end{lem}

\begin{proof}[Proof of Theorem \ref{thmH1}]

We define the metric $h(k)$ on $E$ by
\begin{equation}\label{eqH9}h(k)=\sum_{j=0}^{m}
k^{j-m}\widetilde{g}_{j}=(\sum_{j=0}^{m}
k^{j-m}\Psi_{j})h.\end{equation} Let $B_{k}(h(k),\omega)$ be the
Bergman kernel of $E\otimes L^k$ with respect to the $L^2$-metric
defined by the hermitian metric $h(k)\otimes \sigma^k$ on
$E\otimes L^k$ and the volume form $\frac{ \omega^m}{m!}$ on $X$.
Therefore, if $s_{1},..., s_{N}$ is an orthonormal basis for
$H^0(X,E\otimes L^k)$ with respect to the $L^2(H(k)\otimes
\sigma^k,\frac{ \omega^m}{m!})$, then
\begin{equation}\label{eqH10}B_{k}(h(k),\omega)=\sum s_{i}\otimes s_{i}^{*_{h(k)\otimes
\sigma^k}},\end{equation} We define $\widetilde{B}_{k}(h,\omega)$
as follow
\begin{equation}\label{eqH11}\widetilde{B}_{k}(h,\omega)=\sum s_{i}\otimes s_{i}^{*_{h\otimes
\sigma^k}}.\end{equation} Let
$\widehat{s_{1}},....,\widehat{s_{N}} $ be the corresponding basis
for $H^{0}(\mathbb{P}E^*,\mathcal{O}_{\mathbb{P}E^*}(1)\otimes
L^{k})$. Hence,
$$\int_{\mathbb{P}E^*} \langle \widehat{s_{i}}, \widehat{s_{j}}
\rangle_{g \otimes \sigma^k}d\mu_{g,k}=\int_{\mathbb{P}E^*}
\langle \widehat{s_{i}}, \widehat{s_{j}} \rangle_{g \otimes
\sigma^k}(\sum_{j=0}^{m} k^j f_{j})d\mu_{g}$$
$$= \int_{\mathbb{P}E^*} \langle \widehat{s_{i}}, \widehat{s_{j}}
\rangle_{g(k) \otimes \sigma^k}d\mu_{g}=C_{r} \int_{X}\langle
s_{i}, s_{j} \rangle_{h(k) \otimes
\sigma^k}\frac{\omega^m}{m!}=C_{r}\delta_{ij}.$$ Therefore
$\frac{1}{\sqrt{C_{r}}}\widehat{s_{1}},....,\frac{1}{\sqrt{C_{r}}}\widehat{s_{N}}
$ is an orthonormal basis for
$H^{0}(\mathbb{P}E^*,\mathcal{O}_{\mathbb{P}E^*}(1)\otimes L^{k})$
with respect to $L^2(g \otimes \sigma^k, d\mu_{k,g})$. Hence Lemma
\ref{lemH4} implies
$$C_{r}\rho_{k}(g)=Tr \big(  \lambda(v^*,h)
\widetilde{B}_{k}(h,\omega)       \big).$$ Now, in order to
conclude the proof, it suffices to show that there exist smooth
endomorphisms $A_{i} \in \Gamma(X,E)$ such that
$$ \widetilde{B}_{k}(h,\omega)\sim k^m+A_{1}k^{m-1}+...              .$$
Let $B_{k}(h,\omega)$ be the Bergman kernel of $E\otimes L^k$ with
respect to the $L^2(h\otimes \sigma^k)$. A fundamental result on
the asymptotics of the Bergman kernel (\cite{C}, \cite{Z}) states
that there exists an asymptotic expansion
$$B_{k}(h,\omega)\sim k^m+B_{1}(h)k^{m-1}+...,$$ where $$B_{1}(h)=
\frac{i}{2\pi} \Lambda F_{(E,h)}+ \frac{1}{2} S(\omega) I_{E}
.$$(See also \cite{BBS},\cite{W2}.) Moreover this expansion holds
uniformly for any $h$ in a bounded family. Therefore, we can
Taylor expand the coefficients $B_{i}(h)$'s. We conclude that for
endomorphisms $\Phi_{1},...,\Phi_{M}$,
$$B_{k}(h(I+\sum_{i=0}^{M} k^{-i}\Phi_{i}),\omega)\sim
k^m+B_{1}(h)k^{m-1+...}$$Note that $B_{1}(h)$ in the above
expansion does not depend on $\Phi_{i}$'s and is given as before
by $$B_{1}(h)= \frac{i}{2\pi} \Lambda F_{(E,h)}+ \frac{1}{2}
S(\omega) I_{E}  .$$On the other hand
$$B_{k}(h(k),\omega)=\sum s_{i}\otimes
s_{i}^{*_{\widetilde{g(k)}\otimes \sigma^k}}=(\sum s_{i}\otimes
s_{i}^{*_{h\otimes \sigma^k}})(\sum_{j=0}^m k^{j-m}\Psi_{j})$$
$$=\widetilde{B}_{k}(h,\omega)(\sum_{j=0}^m
k^{j-m}\Psi_{j}).$$ Therefore,
$$\widetilde{B}_{k}(h,\omega)=B_{k}(h(k),\omega)(\sum_{j=0}^m k^{j-m}\Psi_{j})^{-1} \sim k^m+(B_{1}(h)-\Psi_{m-1})k^{m-1}+...$$ Notice that Proposition \ref{propH4}
implies that if $h$ and $\omega$ vary in a bounded family and
$\omega$ is bounded from below, then $\Psi_{1},..,\Psi_{m}$ vary
in a bounded family. Therefore the asymptotic expansion that we
obtained for $\widetilde{B}_{k}(h,\omega)$ is uniform as long as
$h$ and $\omega$ vary in a bounded family and $\omega$ is bounded
from below.

\end{proof}

\begin{prop}\label{propH3}

Suppose that  $\omega_{\infty} \in 2\pi c_{1}(L)$ be a K\"ahler
form with constant scalar curvature and $h_{\textrm{HE}}$ be a
Hermitian-Einstein  metric on $E$, i.e.
$$\Lambda_{\omega_{\infty}} F_{(E,h_{\textrm{HE}})}=\mu I_{E},$$ where
$\mu$ is the $\omega_{\infty}-$slope of the bundle $E$. We have

\begin{align*}A_{1,1}&:=\frac{d}{dt}\Big|_{t=0}A_{1}(h_{\textrm{HE}}(I+t\phi),\omega_{\infty}+it\overline{\partial}\partial
\eta)\\&=\frac{r+1}{2r}\mathcal{D}^*\mathcal{D}\eta I_{E} +
\frac{i}{2\pi}\big((\Lambda_{\omega_{\infty}}\overline{\partial}\partial
\Phi+\Lambda^2_{\omega_{\infty}}(F_{h_{\textrm{HE}}}\wedge
(i\overline{\partial}\partial \eta)
))\\&\qquad-\frac{1}{r}tr(\Lambda_{\omega_{\infty}}\overline{\partial}\partial
\Phi)+\Lambda^2_{\omega_{\infty}}(F_{h_{\textrm{HE}}}\wedge
(i\overline{\partial}\partial \eta) )\big),\end{align*} where
$\mathcal{D}^*\mathcal{D}$ is Lichnerowicz operator (cf.
\cite[Page 515]{D3}).
\end{prop}

\begin{proof}

Define
$$f(t)=\Lambda_{\omega_{\infty}+it\overline{\partial}\partial
\eta}F_{(h_{\textrm{HE}}(I+t\phi))}$$ Therefore, we have
$$mF_{(h_{\textrm{HE}}(I+t\phi))} \wedge
(\omega_{\infty}+it\overline{\partial}\partial \eta)^{m-1}= f(t)
(\omega_{\infty}+it\overline{\partial}\partial \eta)^{m}.$$
Differentiating with respect to $t$ at $t=0$, we obtain

$$m \overline{\partial}\partial
\Phi \wedge \omega_{\infty}^{m-1}+m(m-1)F_{h_{\textrm{HE}}}\wedge
(i\overline{\partial}\partial \eta)\wedge
\omega_{\infty}^{m-2}=f'(0)\omega_{\infty}^m+mf(0)(i\overline{\partial}\partial
\eta)\wedge \omega_{\infty}^{m-1}.$$ Since $f(0)=\mu I_{E}$, we
get
$$f'(0)=\Lambda_{\omega_{\infty}}\overline{\partial}\partial \Phi
+\Lambda^2_{\omega_{\infty}}(F_{h_{\textrm{HE}}}\wedge
(i\overline{\partial}\partial \eta))-\mu
\Lambda_{\omega_{\infty}}(i\overline{\partial}\partial \eta)I_{E}.
$$
On the other hand (cf. \cite[pp. 515, 516]{D3}.)
$$\frac{d}{dt}\Big|_{t=0}S(\omega_{\infty}+it\overline{\partial}\partial
\eta)=\mathcal{D}^*\mathcal{D}\eta.$$

\end{proof}

\begin{cor}\label{corH2}
Suppose that $Aut(X,L)/\mathbb{C}^*$ is discrete and $E$ is
stable. Then the map $A_{1,1}: \Gamma_{0}(End(E)) \rightarrow
\Gamma_{0}(End(E)) $ is an isomorphism, where $\Gamma_{0}(End(E))$
is the space of smooth endomorphisms $\Phi \in E$ such that
$\int_{X}tr(\Phi)\omega_{\infty}^m=0$.

\end{cor}

\begin{proof}
First, notice that
$\Gamma_{0}(End(E))=\Gamma_{00}(End(E))\bigoplus
C^{\infty}_{0}(X)$, where $\Gamma_{00}(End(E))$ is the space of
trace-free endomorphisms of $E$ and $C^{\infty}_{0}(X)$ is  the
space of smooth functions $\eta$ on $X$ such that
$\int_{X}\eta\omega_{\infty}^m=0$.
Assume that $A_{1,1}(\Phi,
\eta)=0$, where $\Phi \in \Gamma_{00}(End(E))$ and $\eta \in
C^{\infty}_{0}(X)$. Hence

$$\frac{r+1}{2r}\mathcal{D}^*\mathcal{D}\eta =0, \,\,\, \textrm{ and}$$
$$\frac{i}{2\pi}\big((\Lambda_{\omega_{\infty}}\overline{\partial}\partial
\Phi+\Lambda^2_{\omega_{\infty}}(F_{h_{\textrm{HE}}}\wedge
(i\overline{\partial}\partial \eta)
))-\frac{1}{r}tr(\Lambda_{\omega_{\infty}}\overline{\partial}\partial
\Phi)+\Lambda^2_{\omega_{\infty}}(F_{h_{\textrm{HE}}}\wedge
(i\overline{\partial}\partial \eta) )\big)=0$$ Since
$Aut(X,L)/\mathbb{C}^*$ is discrete, the first equation implies
that $\eta $ is constant and therefore $\eta=0$. So, the second
equation reduces to the following
$$\Lambda_{\omega_{\infty}}\overline{\partial}\partial \Phi
-\frac{1}{r}tr(\Lambda_{\omega_{\infty}}\overline{\partial}\partial
\Phi)=0$$ It implies that
$$\Lambda_{\omega_{\infty}}\overline{\partial}\partial \Phi=0,$$
since $\Phi$ is traceless. Hence simplicity of $E$ implies that
$\Phi=0$ (cf. \cite{K}).

In order to prove surjectivity let $\Psi \in \Gamma_{0}(End(E)).$
We know that the map $$\eta \in C^{\infty}_{0}\rightarrow
\mathcal{D}^*\mathcal{D}\eta \in C^{\infty}_{0}$$is surjective
since $Aut(X,L)/\mathbb{C}^*$ is discrete (cf. \cite[pp. 515,
516]{D3}). Hence we can find $\eta_{0} $ such that
$\mathcal{D}^*\mathcal{D}\eta_{0}=tr(\Psi).$ On the other hand
$$\frac{i}{2\pi}\big(\Lambda^2_{\omega_{\infty}}(F_{h_{\infty}}\wedge (i\overline{\partial}\partial \eta_{0})
)-\frac{1}{r}tr(\Lambda^2_{\omega_{\infty}}(F_{h_{\infty}}\wedge
(i\overline{\partial}\partial \eta_{0})
)\big)+\Psi-\frac{1}{r}tr(\Psi) \in \Gamma_{0}(End(E)).$$ The map
$$\Phi \in \Gamma_{0}(End(E))\rightarrow
\frac{i}{2\pi}\Lambda_{\omega_{\infty}}\overline{\partial}\partial
\Phi \in \Gamma_{0}(End(E))$$ is surjective since $E$ is simple.
Hence, we can find $\phi_{0}$ such that
$A_{1,1}(\phi_{0},\eta_{0})=\Psi.$

\end{proof}

\section{Constructing Almost Balanced Metrics}
Let $h_{\infty}$ be a hermitian metric on $L$ such that
$\omega_{\infty}=Ric(h_{\infty}) $ be a K\"ahler form with
constant scalar curvature and $h_{\textrm{HE}}$ be the
corresponding Hermitian-Einstein metric on $E$, i.e.
$$\Lambda_{\omega_{\infty}} F_{(E,h_{\textrm{HE}})}=\mu I_{E},$$ where
$\mu$ is the slope of the bundle $E$. Let
$\omega_{0}=Ric(\widehat{h_{\textrm{HE}}})$. After tensoring by
high power of $L$, we can assume without loss of generality that
$\omega_{0}$ is a K\"ahler form on $\mathbb{P}E^*$. We fix an
integer $a \geq 4.$ In order to prove the following, we use ideas
introduced by Donaldson in (\cite[Theorem 26]{D3})

\begin{thm}\label{thmH3}

Suppose $Aut(X,L)$ is discrete. There exist smooth functions
$\eta_{1},\eta_{2},...$ on $X$ and smooth endomorphisms
$\Phi_{1},\Phi_{2},...$ of $E$ such that for any positive integer
$q$ if $$\nu_{k,q}=\omega_{\infty}+i\overline{\partial}\partial
(\sum_{j=1}^q k^{-j}\eta_{j})$$ and
$$h_{k,q}=h_{\textrm{HE}}(I_{E}+\sum_{j=1}^q k^{-j}\Phi_{j}),$$ then

\begin{equation}\label{eqH14}\widetilde{B}_{k}(h_{k,q},\nu_{k,q})=
\frac{C_{r}N_{k}}{k^{-m}V_{k}}(I_{E}+\delta_{q}),\end{equation}
where $||\delta_{q}||_{C^{a+2}}=O(k^{-q-1})$ and
$V_{k}=Vol(\mathbb{P}E^*,\mathcal{O}_{\mathbb{P}E^*}(1)\otimes
L^{k})$ is a topological invariant.

\end{thm}

\begin{proof}
The error term in the asymptotic expansion is uniformly bounded in
$C^{a+2}$ for all $h$ and $\omega$ in a bounded family. Therefore
there exists a positive integer $s$ depends only on $p$ and $q$
such that
\begin{align}\label{eqH15}A_{p}(h(1+\Phi),\omega+i\overline{\partial}\partial
\eta)&=A_{p}(h,\omega)+\sum_{j=1}^{q} A_{p,j}(\Phi,\eta)\\&
\qquad+O(||(\Phi,\eta)||_{C^s}^{q+1}),\notag\end{align} where $
A_{p,j}$ are homogeneous polynomials of degree $j$ , depending on
$h$ and $\omega$, in $\Phi$ and $\eta$ and its covariant
derivatives. Let $\Phi_{1},...,\Phi_{q} $ be smooth endomorphisms
of $E$ and $\eta_{1},...,\eta_{q} $ be smooth functions on $X$. We
have \begin{align}\label{eqH16}&A_{p}(h(1+\sum_{j=1}^q
k^{-j}\Phi_{j}),\omega+i\overline{\partial}\partial (\sum_{j=1}^q
k^{-j}\eta_{j}))\\&\qquad
\qquad=A_{p}(h,\omega)+\sum_{j=1}^{q}b_{p,j}k^{-j}+O(k^{-q-1}),\notag\end{align}
where $b_{p,j}$'s are multi linear expression on $\Phi_{i}$'s and
$\eta_{i}$'s.

Hence
\begin{align}\label{eqH17}&\widetilde{B}_{k}(h(1+\sum_{j=1}^q
k^{-j}\Phi_{j}),\omega+i\overline{\partial}\partial (\sum_{j=1}^q
k^{-j}\eta_{j}))\\&=k^m+A_{1}(h,\omega)k^{m-1}+....\notag\\&\qquad+(A_{q}(h,\omega)+b_{q-1,1}+...
+b_{1,q-1})k^{m-q}+O(k^{m-q-1}).\notag\end{align} We need to
choose $\Phi_{j}$ and $\eta_{j}$ such that coefficients of
$k^m,...k^{m-q}$ in the right hand side of \eqref{eqH17} are
constant. Donaldson's key observation is that $\eta_{p}$ and
$\phi_{p}$ only appear in the coefficient of $k^{m-p}$ in the form
of $A_{1,1}(\phi_{p}, \eta_{p})$. Hence, we can do this
inductively. Assume that we choose
$\eta_{1},\eta_{2},...\eta_{p-1}$  and
$\Phi_{1},\Phi_{2},...,\Phi_{p-1}$ so that  the coefficients of
$k^m,...k^{m-p+1}$ are constant. Now we need to choose $\eta_{p}$
and $\Phi_{p}$ such that the coefficient of $k^{m-p}$ is constant.
This means that  we need to solve the equation
\begin{equation}\label{eqH12}A_{1,1}(\Phi_{p},\eta_{p})-c_{p}I_{E}= P_{p-1},\end{equation}
for $\Phi_{p},\eta_{p}$ and the constant $c_{p}$. In this equation
$P_{p-1}$ is determined by $\Phi_{1},...,\Phi_{p-1}$ and
$\eta_{1},...,\eta_{p-1}$. Corollary \ref{corH2} implies that we
can always solve the equation \eqref{eqH12}.

\end{proof}

\begin{cor}\label{corH3}

For any positive integer $q$, there exist hermitian metrics
$g_{k,q}$ on $\mathcal{O}_{\mathbb{P}E^*}(1)$ and K\"ahler forms
$\nu_{k,q}$ on $X$ in the class of $2\pi c_{1}(L)$ so that
$$\rho_{k}(g_{k,q},\nu_{k,q})=\frac{N_{k}}{k^{-m}V_{k}}(1+\epsilon_{k,q}),$$
where $||\epsilon_{k,q}||_{C^{a+2}}=O(k^{-q-1}).$ Moreover,
\begin{equation}\label{neq6}||\omega_{g_{k,q}}+k\nu_{k,q}-(\omega_{0}+k\omega_{\infty})||_{C^{a}(\omega_{0}
+k\omega_{\infty})}=O(k^{-1}).\end{equation}

\end{cor}

\begin{proof}
Let $g_{k,q}=\widehat{h_{k,q}}.$ By Theorem \ref{thmH3}, We have
\begin{align*}\rho_{k}(g_{k,q},\nu_{k,q})&=\frac{N_{k}}{k^{-m}V_{k}}Tr(\lambda(v^*,h_{k,q})(I_{E}+\delta_{q}))\\&=\frac{N_{k}}{k^{-m}V_{k}}(1+Tr(\lambda(v^*,h_{k,q})\delta_{q}))).\end{align*}
The first part of corollary is proved, since $h_{k,q}$ is bounded
and $||\delta_{k,q}||_{C^{a+2}}=O(k^{-q-1})$. Define
$\widetilde{\omega_{0}}=\omega_{0}+k\omega_{\infty}$. For the
second part, we have
\begin{align*}||\omega_{g_{k,q}}+k\nu_{k,q}-(\omega_{0}+k\omega_{\infty})||_{C^{a}(\widetilde{\omega_{0}}
)} &\leq
||\omega_{g_{k,q}}-\omega_{0}||_{C^{a}(\widetilde{\omega_{0}}
)}+k||\nu_{k,q}-\omega_{\infty}||_{C^{a}(\widetilde{\omega_{0}}
)}\\&\leq ||\omega_{g_{k,q}}-\omega_{0}||_{C^{a}(\omega_{0}
)}+k||\nu_{k,q}-\omega_{\infty}||_{C^{a}(k\omega_{\infty}
)}\\&=||\omega_{g_{k,q}}-\omega_{0}||_{C^{a}(\omega_{0}
)}+||\nu_{k,q}-\omega_{\infty}||_{C^{a}(\omega_{\infty}
)}\\&=O(k^{-1}).\end{align*} Notice that by definition, we have
\begin{align*}&||\omega_{g_{k,q}}-\omega_{0}||_{C^a(\omega_{0})}=O(k^{-1}),\\&||\nu_{k,q}-\omega_{\infty}||_{C^a(\omega_{\infty})}=O(k^{-1}).\end{align*}

\end{proof}

\section{Proof of the main theorem}

In this section, we prove Theorem \ref{thm2}. In order to do that,
we want to apply Theorem \ref{newthm2}. Hence, we need to
construct a sequence of almost balanced metrics on
$\mathbb{P}E^*,\mathcal{O}_{\mathbb{P}E^*}(1)\otimes L^{\otimes
k}$. Also, we need to show that $\mathbb{P}E^*$ has no nontrivial
holomorphic vector fields.

\begin{prop}\label{newprop2}
Let $E $ be a holomorphic vector bundle over a compact K\"ahler
manifold $X$. Suppose that $X$ has no nonzero holomorphic vector
fields. If $E$ is stable, then  $\mathbb{P}E^*$ has no nontrivial
holomorphic vector fields.

\end{prop}

\begin{proof}

Let $TF$ be the sheaf of tangent vectors to the fibre of $\pi$. We
have the following exact sequence over $\mathbb{P}E$:

$$0\rightarrow TF \rightarrow T \mathbb{P}E^* \rightarrow \pi^* TX \rightarrow 0.$$
This gives the long exact sequence
$$0\rightarrow H^0(\mathbb{P}E^* ,TF) \rightarrow H^0(\mathbb{P}E^* ,T \mathbb{P}E^*) \rightarrow H^0(\mathbb{P}E^* ,\pi^* TX ) \rightarrow \dots$$
Since $H^0(\mathbb{P}E^* ,\pi^* TX )=0$ , we have
$$H^0(\mathbb{P}E^* ,TF) \simeq H^0(\mathbb{P}E^* ,T \mathbb{P}E^*)$$ On
the other hand, $\pi_{*}TF $ may be identified with the sheaf of
trace free endomorphisms of $E$. Therefore by simplicity of $E$
(cf. \cite{K})
$$H^0(\mathbb{P}E^*,TF) \simeq H^0(X, \pi_{*} TF) =0.$$

\end{proof}

\begin{proof}[Proof of Theorem \ref{thm2}]

Since Chow stability is equivalent to the existence of balanced
metric, it suffices to show that
$(\mathbb{P}E^*,\mathcal{O}_{\mathbb{P}E^*}(1)\otimes \pi^*L^{k})$
admits balanced metric for $k \gg 0.$ Fix a positive integer $q$.
From now on we drop all indexes $q$ for simplicity. Let
$\sigma_{k}=\sigma_{k,q}$ be a metric on $L$ such that
$Ric(\sigma_{k})=\nu_{k}$, where $\nu_{k}=\nu_{k,q}$ is the one in
the statement of Theorem \ref{thm3}. Let $t_{1},...,t_{N}$ be an
orthonormal basis for
$H^{0}(\mathbb{P}E^*,\mathcal{O}_{\mathbb{P}E^*}(1)\otimes L^{k})$
w.r.t. $L^2(g_{k}\otimes \sigma_{k}^{\otimes k},
\frac{(\omega_{g_{k}}+k\nu_{k})^{m+r-1}}{(m+r-1)!}).$ Thus,
Corollary \ref{corH3} implies $$\sum |t_{i}|^2_{g_{k}\otimes
\sigma_{k}^{\otimes k}}=\frac{N_{k}}{V_{k}}(1+\epsilon_{k}).$$
Define $g_{k}^{'}=\frac{V_{k}}{N_{k}}(1+\epsilon_{k})^{-1}g_{k}$.
We have $$\sum |t_{i}|^2_{g_{k}^{'}\otimes \sigma_{k}^{\otimes
k}}=1.$$ This implies that the metric $g_{k}^{'}$ is the
Fubini-Study metric on $\mathcal{O}_{\mathbb{P}E^*}(1)\otimes
L^{k}$ induced by the embedding $\iota_{\underline{t}}:
\mathbb{P}E^* \rightarrow \mathbb{P}^{N-1},$ where
$\underline{t}=(t_{1},...,t_{N})$. We prove that this sequence of
embedding is almost balanced of order $q$, i.e
$$\int_{\mathbb{P}E^*} \langle t_{i},
t_{j}\rangle_{g_{k}^{'}\otimes \sigma_{k}^{\otimes
k}}\frac{(\omega_{g_{k}^{'}}+k\nu_{k})^{m+r-1}}{(m+r-1)!}=D^{(k)}\delta_{ij}+M_{ij},$$
where $M^{(k)}=[M_{ij}]$ is a trace free hermitian matrix,
$D^{(k)} \rightarrow C_{r}$ as $k \rightarrow \infty$ and
$||M^{(k)}||_{\textrm{op}}=O(k^{-q-1}).$

\begin{align*}M_{ij}^{(k)}&=\int_{\mathbb{P}E^*} \langle
t_{i},t_{j}\rangle_{g_{k}^{'}\otimes\sigma_{k}^{\otimes
k}}\frac{(\omega_{g_{k}^{'}}+k\nu_{k})^{m+r-1}}{(m+r-1)!}\\&\qquad-\frac{V_{k}}{N_{k}}\int_{\mathbb{P}E^*}\langle
t_{i}, t_{j}\rangle_{g_{k}\otimes \sigma_{k}^{\otimes
k}}\frac{(\omega_{g_{k}}+k\nu_{k})^{m+r-1}}{(m+r-1)!}\\&=\frac{V_{k}}{N_{k}}
\int_{\mathbb{P}E^*}\langle t_{i}, t_{j}\rangle_{g_{k}\otimes
\sigma_{k}^{\otimes
k}}(f_{k}(1+\epsilon_{k})^{-1}-1)\frac{(\omega_{g_{k}}+k\nu_{k})^{m+r-1}}{(m+r-1)!},\end{align*}
where
$$(\omega_{g_{k}^{'}}+k\nu_{k})^{m+r-1}=f_{k}(\omega_{g_{k}}+k\nu_{k})^{m+r-1}.$$
By a unitary change of basis, we may assume without loss of
generality that the matrix $M^{(k)}$ is diagonal. Thus
$$||M^{(k)}||_{\textrm{op}} \leq
\frac{V_{k}}{N_{k}}||f_{k}(1+\epsilon_{k})^{-1}-1||_{L^{\infty}}.$$
On the other hand,
\begin{align*}||\omega_{g_{k}^{'}}-\omega_{g_{k}}||_{C^0(\omega_{0})}&=||\overline{\partial}\partial
\log(1+\epsilon_{k})||_{C^0(\omega_{0})}\\&\leq
||\log(1+\epsilon_{k})||_{C^2(\omega_{0})}.\\&\leq
-\log(1-C||\epsilon_{k}||_{C^2(\omega_{0})})\\&=
O(k^{-q-1}).\end{align*} Therefore,
\begin{align*}||f_{k}-1||_{\infty}&=\Big|\frac{\omega^{m+r-1}_{g_{k}^{'}}-\omega_{g_{k}}^{m+r-1}}{\omega_{g_{k}}^{m+r-1}}\Big|=
\Big|\frac{\omega^{m+r-1}_{g_{k}^{'}}-\omega_{g_{k}}^{m+r-1}}{\omega_{0}^{m+r-1}}\frac{\omega_{0}^{m+r-1}}{\omega_{g_{k}}^{m+r-1}}\Big|\\&\leq
Ck^{-q-1}\Big|\frac{\omega_{0}^{m+r-1}}{\omega_{g_{k}}^{m+r-1}}\Big|.
\end{align*}
This implies  that $$||f_{k}-1||_{\infty}\leq Ck^{-q-1},$$ since
$\Big|\frac{\omega_{0}^{m+r-1}}{\omega_{g_{k}}^{m+r-1}}\Big|$ is
bounded. Hence
$$||f_{k}(1+\epsilon_{k})^{-1}-1||\leq C^{'}k^{-q-1}.$$ Therefore
$$||M^{(k)}||_{\textrm{op}} =O(k^{-q-1}).$$

Proposition \ref{newprop2} implies that $\mathbb{P}E^*$ has no
nontrivial holomorphic vector fields. Therefore, applying Theorem
\ref{newthm2} and \eqref{neq6} conclude the proof.

\end{proof}

\begin{rem}
Since $\mathbb{P}E^*$ has no nontrivial holomorphic vector fields,
the sequence $\omega_{k}^{\textrm{bal}}$ of balanced metrics in
the class of $\mathcal{O}_{\mathbb{P}E^*}(1)\otimes L^{\otimes k}$
is unique. Define the sequence of metrics
$\Omega_{k}=\omega_{k}^{\textrm{bal}}-k\omega_{\infty}$ in the
class of $\mathcal{O}_{\mathbb{P}E^*}(1).$ A natural question is
whether $\Omega_{k}$ converges to $\omega_{0}$ as $k \rightarrow
\infty$. If $\dim_{\mathbb{C}}X=1$, then it is easy to see that
$\Omega_{k}$ converges to $\omega_{0}$ in $C^{\infty}$-norm as $k
\rightarrow \infty$. In general, we have
\begin{align*}||\Omega_{k}-\omega_{0}||_{C^{p}(\omega_{0})}&=||\omega_{k}^{\textrm{bal}}-\widetilde{\omega_{0}}||_{C^{p}(\omega_{0})}\\&\leq||\omega_{k}^{\textrm{bal}}-\widetilde{\omega_{k}}||_{C^{p}(\omega_{0})}+||\widetilde{\omega_{k}}-\widetilde{\omega_{0}}||_{C^{p}(\omega_{0})}\end{align*}
The first term has order of $O(k^{-1})$ by \eqref{eq15} and the
second term is of order $O(1)$. Therefore, in higher dimension one
gets
$$||\Omega_{k}-\omega_{0}||_{C^{p}(\omega_{0})}=O(1) ,$$
for any positive integer $p$. It is not clear that whether one can
find a sharper estimate with these methods.

\end{rem}


\end{document}